\documentclass{amsart}
\usepackage{amscd}
\usepackage{amsmath}
\usepackage{amssymb}
\usepackage{amsthm}
\usepackage{latexsym}
\usepackage{pstricks}

\newtheorem{theorem}{Theorem}[section]
\newtheorem{lemma}[theorem]{Lemma}
\newtheorem{proposition}[theorem]{Proposition}
\newtheorem{corollary}[theorem]{Corollary}

\theoremstyle{definition}
\newtheorem{definition}[theorem]{Definition}

\newtheorem*{acknowledgement}{Acknowledgment}
\theoremstyle{remark}
\newtheorem{remark}[theorem]{Remark}
\numberwithin{equation}{section}

\DeclareMathOperator{\pd}{pd}

\DeclareMathOperator{\height}{height}

\DeclareMathOperator{\indeg}{indeg}

\DeclareMathOperator{\reg}{reg}
\DeclareMathOperator{\arithdeg}{arithdeg}
\DeclareMathOperator{\ara}{ara}
\DeclareMathOperator{\chara}{char}

\DeclareMathOperator{\cone}{cone}

\newcommand{\calH}{{\mathcal{H}}}

\def\pot#1#2{#1[\kern-0.28ex[#2]\kern-0.28ex]}

\def\Implies{\ifmmode\Longrightarrow \else
     \unskip${}\Longrightarrow{}$\ignorespaces\fi}
\def\implies{\ifmmode\Rightarrow \else
     \unskip${}\Rightarrow{}$\ignorespaces\fi}
\def\iff{\ifmmode\Longleftrightarrow \else
     \unskip${}\Longleftrightarrow{}$\ignorespaces\fi}

\let\:=\colon
\def\opn#1#2{\def#1{\operatorname{#2}}} 

\let\epsilon\varepsilon
\let\kappa=\varkappa
\opn\Gin{Gin}

\opn\inii{in} \opn\inim{inm} \opn\rate{rate}

\begin{document}
\title[Arithmetical rank of squarefree monomial ideals ]
{Arithmetical rank of squarefree monomial ideals 
generated by five elements or with arithmetic degree four }
\author[K. Kimura]{Kyouko Kimura}
\address[Kyouko Kimura]{Department of Mathematics, 
         Faculty of Science, 
         Shizuoka University, 836 Ohya, Suruga-ku, Shizuoka 422--8529, Japan}
\email{skkimur@ipc.shizuoka.ac.jp}
\author[G. Rinaldo]{Giancarlo Rinaldo}
\address[Giancarlo Rinaldo]{Dipartimento di Matematica, 
Universita' di Messina, Salita Sperone, 31. S. Agata, 
Messina 98166, Italy}
\email{rinaldo@dipmat.unime.it}
\author[N. Terai]{Naoki Terai}
\address[Naoki Terai]{Department of Mathematics, 
         Faculty of Culture and Education, 
         Saga University, Saga 840--8502, Japan}
\email{terai@cc.saga-u.ac.jp}
\begin{abstract}
Let $I$ be a squarefree monomial ideal of a polynomial ring $S$. 
In this paper, we prove that the arithmetical rank of $I$ is equal to 
the projective dimension of $S/I$ when one of the following conditions 
is satisfied: (1) $\mu (I) \leq 5$; (2) $\arithdeg I \leq 4$. 
\end{abstract}
\maketitle
\section{Introduction}
Let $S$ be a polynomial ring over a field $K$ and 
$I$ a squarefree monomial ideal of $S$. 
We denote by $G(I)$, the minimal set of monomial generators of $I$. 
The \textit{arithmetical rank} of $I$, denoted by $\ara I$, is defined by 
the minimal number $r$ of elements $a_1, \ldots, a_r \in S$ such that 
\begin{displaymath}
  \sqrt{(a_1, \ldots, a_r)} = \sqrt{I}. 
\end{displaymath}
When the above equality holds, 
we say that \textit{$a_1, \ldots, a_r$ generate $I$ up to radical}. 
By definition, $\ara I \leq \mu (I)$ holds, 
where $\mu (I)$ denotes the cardinality of $G(I)$. 
On the other hand, Lyubeznik \cite{Ly83} proved that 
\begin{equation}
  \label{eq:ara>pd}
  \pd_S S/I \leq \ara I, 
\end{equation}
where $\pd_S S/I$ denotes the projective dimension of $S/I$ over $S$. 
Since $\height I \leq \pd_S S/I$ always holds, we have 
\begin{displaymath}
  \height I \leq \pd_S S/I \leq \ara I \leq \mu (I). 
\end{displaymath}
Then it is natural to ask when $\ara I = \pd_S S/I$ holds. 
Many authors investigated this problem; see 
\cite{Barile96, Barile08-1, Barile08-2, 
BKMY, BariTera08, BariTera09, Kim_h2CM,  KTYdev1, KTYdev2, Ku, Mo,  SchmVo}. 
In particular in \cite{KTYdev1, KTYdev2}, 
it was proved that the equality $\ara I = \pd_S S/I$ holds for the two cases that
$\mu (I)-\height I \le 2$ and 
$\arithdeg I-\indeg I \le 1$. 
Here $\arithdeg I$ denotes the \textit{arithmetic degree} of $I$, 
which is equal to the number of minimal primes of $I$, 
and $\indeg I$ denotes the \textit{initial degree} of $I$; see Section 2. 
As a  result we know that  $\ara I = \pd_S S/I$ for squarefree monomial ideals $I$
with $\mu (I) \le 4 $  or with $\arithdeg I \le 3 $. 

\par
In this paper, hence, we concentrate our attention on the following two cases:
$\mu (I) = 5 $ and  $\arithdeg I = 4$. 
And the main result of this paper is as follows: 
\begin{theorem}
  \label{claim:MainResult}
  Let $I$ be a squarefree monomial ideal of $S$. 
  Suppose that $I$ satisfies one of the following conditions$:$ 
  \begin{enumerate}
  \item $\mu (I) \leq 5$. 
  \item $\arithdeg I \leq 4$. 
  \end{enumerate}
  Then 
  \begin{displaymath}
    \ara I = \pd_S S/I. 
  \end{displaymath}
\end{theorem}
Note that there exists an ideal $I$ with $\mu (I) = 6$ such that 
$\ara I > \pd_S S/I$ when $\chara K \neq 2$; see \cite[Section 6]{KTYdev2}. 

\par
After we recall some definitions and properties of 
Stanley--Reisner ideals and hypergraphs in Sections 2 and 3, 
we give a combinatorial characterization
for a squarefree monomial ideal $I$ with  $\pd_S  S/I =\mu (I)-1$
using hypergraphs in Section 4.
It is necessary because we must use the fact that projective dimension of $S/I$ is characteristic-free 
for a  squarefree monomial ideal $I$ with $\mu(I)=5$ in Section 6.

\par
In the case that $\mu (I) \le 4 $
all the squarefree monomial ideals are classified using hypergraphs in \cite{KTYdev1, KTYdev2}.
But in the case that $\mu (I) =5$, $\height I=2$ and $\pd S/I=3$, which is an essential difficult part
for our problem,
a similar classification is practically  impossible because of their huge number.
According to a computer, there are about $2.3 \cdot 10^6$ hypergraphs 
corresponding to such ideals. 
We need a reduction.
Hence we focus on the set of the most ``general" members among them.
We call it a {\it generic set}.
In Section 5
we give a formal definition of a generic set and prove that
it is enough to show $\ara I = \pd_S S/I$ for each member $I$ of the generic set.
But we cannot obtain the generic set in our case without a computer.
In Section 6, we give an algorithm to find a generic set 
for the ``connected" squarefree monomial ideals $I$ with
$\mu (I) =5$, $\height I=2$ and $\pd S/I=3$ using {\it CoCoA} and {\it Nauty}.
As a result of computation, we found that the generic set consists of just three ideals. 
In Section 7, we prove $\ara I = \pd_S S/I$ when $\mu (I) =5$
by showing the same equality for all three members of the generic set.

\par
Finally in Section 8, we focus on the squarefree monomial ideals $I$ with $\arithdeg I \leq 4$. 
Here we use another reduction.
In terms  of simplicial complexes as shown in \cite{BariTera08, Kim_h2CM},
we may remove 
a face with a free vertex from a simplicial complex $\Delta$
when we consider the problem whether $\ara I_{\Delta} = \pd_S S/I_{\Delta}$ holds. 
We translate it in terms of hypergraphs and after such reduction 
we show  $\ara I = \pd_S S/I$ for remaining ideals $I$ 
with  $\arithdeg I \leq 4$.

\section{Preliminaries}
In this section, we recall some definitions and properties which 
are needed to prove Theorem \ref{claim:MainResult}. 

\par
Let $M$ be a Noetherian graded $S$-module and 
\begin{displaymath}
  F_{\bullet} : 0 \longrightarrow \bigoplus_{j \geq 0} S(-j)^{{\beta}_{p, j}} 
    \longrightarrow \cdots \longrightarrow 
    \bigoplus_{j \geq 0} S(-j)^{{\beta}_{0, j}} \longrightarrow M 
    \longrightarrow 0
\end{displaymath}
a graded minimal free resolution of $M$ over $S$, 
where $S(-j)$ is a graded free $S$-module whose $k$th piece is 
given by $S_{k-j}$. 
Then ${\beta}_{i, j} = {\beta}_{i, j} (M)$ is called 
a \textit{graded Betti number} of $M$ 
and ${\beta}_i = \sum_{j} {\beta}_{i, j}$ is called 
the $i$th \textit{$($total$)$ Betti number} of $M$. 
The projective dimension of $M$ over $S$ is defined  by $p$ and denoted by $\pd_S S/I$ or by $\pd S/I$.
The \textit{initial degree} of $M$ and the \textit{regularity} of $M$ 
are defined by 
\begin{displaymath}
  \indeg M = \min \{ j \; : \; {\beta}_{0, j} (M) \neq 0 \}, 
  \quad \reg M = \max \{ j-i \; : \; {\beta}_{i, j} (M) \neq 0 \}, 
\end{displaymath}
respectively. 

\par
Next, we recall some definitions and properties of 
Stanley--Reisner ideals, especially Alexander duality. 

\par
Let $V =[n]:= \{ 1, 2, \ldots, n \}$. A \textit{simplicial complex} $\Delta$ 
on the vertex set $V$ is a collection of subsets of $V$ with the conditions 
(a) $\{ v \} \in \Delta$ for all $v \in V$; 
(b) $F \in \Delta$ and $G \subset F$ imply $G \in \Delta$. 
An element of $V$ is called a \textit{vertex} of $\Delta$ 
and an element of $\Delta$ is called a \textit{face}. 
A maximal face of $\Delta$ is called a \textit{facet} of $\Delta$. 
Let $F$ be a face of $\Delta$. The \textit{dimension} of $F$, 
denoted by $\dim F$, 
is defined by $|F|-1$, where $|F|$ denotes the cardinality of $F$. 
If $\dim F = i$, then $F$ is called an $i$-face. 
The \textit{dimension} of $\Delta$ is defined by 
$\dim \Delta := \max \{ \dim F : F \in \Delta \}$. 
A simplicial complex which consists of all subsets of its vertex set 
is called a \textit{simplex}. 
Let $u$ be a new vertex and $F \subset V$. 
The \textit{cone from $u$ over $F$} is a simplex on the vertex set 
$F \cup \{u \}$; 
see \cite[Definition 1, p.\  3687]{BariTera08}. 
We denote it by $\cone_u F$. 
Then the union $\Delta \cup \cone_u F$ is 
a simplicial complex on the vertex set $V \cup \{ u \}$. 

\par
The \textit{Alexander dual complex} ${\Delta}^{\ast}$ of $\Delta$ 
is defined by 
  ${\Delta}^{\ast} = \{ F \subset V : V \setminus F \notin \Delta \}$. 
If $\dim \Delta < n-2$, then ${\Delta}^{\ast}$ is also a simplicial complex 
on the same vertex set $V$. 

\par
For a simplicial complex $\Delta$ 
on the vertex set $V = [n]$, 
we can associate a squarefree monomial ideal $I_{\Delta}$ of 
$S = K[x_1, \ldots, x_n]$ which is generated by all products 
$x_{i_1} \cdots x_{i_s}$, $1 \leq i_1 < \cdots < i_s \leq n$ with 
$\{ i_1, \ldots, i_s \} \notin \Delta$. 
The ideal $I_{\Delta}$ is called the \textit{Stanley--Reisner ideal} 
of $\Delta$. 
It is well known that the minimal prime decomposition of $I_{\Delta}$ is given by 
\begin{displaymath}
  I_{\Delta} = \bigcap_{F \in \Delta : \text{facet}} P_{F}, 
\end{displaymath}
where $P_G = (x_i : i \in V \setminus G)$ for $G \subset V$. 
On the other hand, for a squarefree monomial ideal $I \subset S$ with 
$\indeg I \geq 2$, there exists a simplicial complex $\Delta$ on the 
vertex set $V = [n]$ such that $I = I_{\Delta}$. 
If $\dim \Delta < n-2$, i.e., $\height I \geq 2$, then we can consider 
the squarefree monomial ideal $I^{\ast} := I_{{\Delta}^{\ast}}$, 
which is called the \textit{Alexander dual ideal} of $I = I_{\Delta}$. 
Then 
\begin{displaymath}
  I^{\ast} = I_{{\Delta}^{\ast}} = (x^{V \setminus F} \; : \; 
    \text{$F \in \Delta$ is a facet}), 
\end{displaymath}
where $x^G = \prod_{i \in G} x_i$ for $G \subset V$. 
It is easy to see that $I^{\ast \ast} = I$, 
$\indeg I^{\ast} = \height I$, and $\arithdeg I^{\ast} = \mu (I)$ hold. 
Moreover, the equality $\reg I^{\ast} = \pd_S S/I$ also holds; 
see \cite[Corollary 1.6]{Terai}. 

\section{Hypergraphs}
For this section, we refer to Kimura, Terai and Yoshida \cite{KTYdev1}, 
\cite{KTYdev2} for more detailed information. 

\par
Let $V = [\mu]$. A \textit{hypergraph} $\calH$ 
on the vertex set $V$ is a collection of subsets of $V$ 
with $\bigcup_{F \in \calH} F = V$. 
The definitions and notations of the vertex, face, and dimension are the 
same as those for a simplicial complex. 
We set $B({\calH}) = \{ v \in V : \{ v \} \in \calH \}$ and $W(\calH)=V \setminus B(\calH)$.
For a hypergraph  $\calH$ on a vertex set $V$,
we define the {\it  $i$-subhypergraph} of $\calH$ by $\calH^{i }=\{ F \in \calH: \  \dim F= i\}$.
We sometimes identify $B(\calH)$ with $\calH^{0}$. 
For $U \subset V(\calH)$, we define the restriction of a hypergraph $\calH$ to $U$ by
$\calH_{U}=\{ F \in \calH: \  F \subset U \}$.

\par
A hypergraph $\calH$ on the vertex set $V$ is called \textit{disconnected} 
if there exist hypergraphs ${\calH}_1, {\calH}_2 \subsetneq \calH$ 
on vertex sets $V_1, V_2 \subsetneq V$, respectively such that 
${\calH}_1 \cup {\calH}_2 = \calH$, $V_1 \cup V_2 = V$, and 
$V_1 \cap V_2 = \emptyset$. 
A hypergraph which is not disconnected is called \textit{connected}. 

\par
Let $I$ be a squarefree monomial ideal of $S=K[x_1, \ldots, x_n]$ 
with $G(I) = \{ m_1, \ldots, m_{\mu} \}$. 
We associate a hypergraph $\calH (I)$ on the vertex set 
$V = [\mu]$ with $I$ by setting 
\begin{displaymath}
  \calH (I) := \big\{ \{ j \in V \; : \; \text{$m_j$ is divisible by $x_i$} \} 
    \; : \; i=1, 2, \ldots, n \big\}. 
\end{displaymath}

\begin{definition}
  \label{defn:def_var}
  Let $F$ be a face of $\calH (I)$. 
  Then there exists a variable $x$ of $S$ such that 
  \begin{equation}
    \label{eq:def_var}
    F = \{ j \in V \; : \; \text{$m_j$ is divisible by $x$} \}. 
  \end{equation}
  We call a variable $x$ of $S$ with condition (\ref{eq:def_var}) 
  a \textit{defining variable} of $F$ (in $\calH (I)$). 
     \par
  Conversely, we say that \textit{a variable $x$ of $S$ defines a face $F$} of 
  $\calH (I)$ if $x$ is a defining variable of $F$. 
\end{definition}

 Note that the choice of a defining variable is not necessarily unique.
 Since the minimal generators of a concrete squarefree monomial ideal $I$ do not have indices,
 we usually regard $\calH (I)$ as a hypergraph with unlabeled vertices.
 If $\calH (I)$ can be regarded as a subhypergraph of
 $\calH (J)$ as unlabeled hypergraphs for two squarefree monomial ideals $I$ and $J$, 
 we write $\calH (I) \subset \calH (J)$ by abuse of language.

\par
On the other hand, we can construct a squarefree monomial ideal from a given 
hypergraph $\calH$ on the vertex set $V = [\mu]$ 
if $\calH$ satisfies the following \textit{separability condition}: 
\begin{displaymath}
  \begin{aligned}
    &\text{For any two vertices $i, j \in V$,} \\
    &\text{there exist faces $F, G \in \calH$ 
    such that $i \in F \setminus G$ and $j \in G \setminus F$. }
  \end{aligned}
\end{displaymath}
The way of construction is: 
first, we assign a squarefree monomial $A_F$ to each face $F \in \calH$ 
such that $A_F$ and $A_G$ are coprime if $F \neq G$. Then we set 
\begin{displaymath}
  I = (\prod_{\genfrac{}{}{0pt}{}{F \in \calH}{j \in F}} A_F \; : \; 
               j = 1, 2, \ldots, \mu), 
\end{displaymath}
which is a squarefree monomial ideal with $\calH (I) = \calH$ 
by virtue of the separability. 
When we assign a variable $x_F$ for each $F \in \calH$, 
we write the corresponding ideal as $I_{\calH}$ in $K[x_F : F \in \calH]$. 

\par
For later use we prove the following proposition:

\begin{proposition}
  \label{claim:hgraph}
  Let $I, I'$ be squarefree monomial ideals of polynomial rings $S, S'$, 
  respectively. 
  Suppose that $\mu (I)=\mu (I')$ and  $\calH (I) \subset \calH (I')$. 
  Then we have $\ara I \leq \ara I'$. 
\end{proposition}
\begin{proof}
  Set $\calH = \calH (I)$, $\calH' = \calH (I')$. 
  We may assume that $I=I_{\calH }$ and  $I'=I_{\calH '}$
  with ${\calH \subset \calH '}$.
  Set $G(I')=\{m'_1, \dots, m'_{\mu}\}$.
  Let $m_i$ for $i=1,2, \dots ,\mu$ be the monomial obtained by
  substitution of 1 to $x_F$ for $F \in \calH ' \setminus \calH $ in $m'_i$.
  We may assume that $G(I)=\{m_1, \dots, m_{\mu}\}$.
  Assume that $q_1', \ldots, q_r'$ generate $I'$ up to radical.
  Let $q_i$ for $i=1,2, \dots ,r$ be the polynomial obtained by
  substitution of 1 to $x_F$ for $F \in \calH ' \setminus \calH $ in $q'_i$.
  We show that $q_1, \ldots, q_r$ generate $I$ up to radical.
  Since $q'_i \in I'$ for $i=1,2, \dots ,r$, we have $q_i \in I$.
  On the other hand, suppose ${m'_i}^p \in ( q_1', \ldots, q_r')$ for some $p \ge 1$.
  Then we have $m_i^p \in ( q_1, \ldots, q_r)$.
  Hence  $q_1, \ldots, q_r$ generate $I$ up to radical.
\end{proof}

\section{Squarefree monomial ideals whose projective dimension is close  to the number of generators}
Let $S= K[x_{1},x_{2},\ldots,x_{n}]$ be the polynomial 
ring in $n$ variables over a field $K$. 
We fix a squarefree monomial ideal $I= (m_{1}, m_{2},  \dots , m_{\mu})$, where 
$G(I)=\{m_{1}, m_{2},  \dots , m_{\mu} \}$ 
is the minimal generating set of monomials for
$I$.
For a squarefree monomial ideal $I$, by the Taylor resolution of $S/I$
we have $\pd_S S/I \le \mu (I)$.
In this section we give a combinatorial characterization
for the squarefree monomial ideal $I$ with  $\pd_S  S/I =\mu (I)-1$
using hypergraphs.

\par
First we consider the condition $\pd_S S/I = \mu (I)$.
Then the following proposition is easy and well known.
\begin{proposition}
\label{claim:mu}
The following conditions are equivalent for a squarefree monomial ideal $I:$ 
\begin{enumerate}
\item[(1)] $\pd_S S/I = \mu (I)$. 
\item[(2)] For the hypergraph $\calH:=\calH(I)$ we have $B(\calH) =V(\calH)$. 
\end{enumerate}
\end{proposition}

\par
By Lyubeznik\cite{Ly83} we have $\pd_S S/I = 
\mbox{\rm cd} I:= \max \{ i: H^{i}_{I}(S) \ne 0 \}$.
Assuming that $\pd_S  S/I \le \mu (I)-1$,
we have  $\pd_S  S/I =\mu (I)-1$ if and only if $H^{\mu -1}_{I}(S) \ne 0$.
We give a combinatorial interpretation for the condition  $H^{\mu -1}_{I}(S) \ne 0$.

\par
Consider the following {\v C}ech complex: 
\begin{eqnarray*}
\lefteqn{
C^{\bullet}
 = 
\bigotimes_{i=1}^{\mu}
(0 
\longrightarrow 
S 
\longrightarrow
S_{m_{i}}
\longrightarrow
0 )
}\\ 
& &  = 
0 
\longrightarrow 
S 
\stackrel{\delta^{1}}{\longrightarrow}
\bigoplus_{1 \le i \le \mu } S_{m_{i}}
\stackrel{\delta^{2}}{\longrightarrow} 
\bigoplus_{1 \le i < j \le \mu } S_{m_{i}m_{j}}
\stackrel{\delta^{3}}{\longrightarrow}
\cdots
\stackrel{\delta^{\mu}}{\longrightarrow}
S_{m_{1}m_{2}\cdots m_{\mu}}
\longrightarrow
0.
\end{eqnarray*}
We describe $\delta^{r+1}$ as follows. 
Put $R:= S_{m_{i_{1}}m_{i_{2}}\cdots m_{i_{r}}}$
and $\{j_{1},j_{2}, \dots , j_{s}\}= \{ 1,2, \cdots ,\mu \}\setminus
\{i_{1},i_{2}, \dots , i_{r}\}$,
where $j_{1}< j_{2} < \dots < j_{s}$ and $r+s=\mu$.
Let $\psi _{j_p}: R \longrightarrow R_{m_{j_p}}$ be a natural map.
For $u \in R$,
we have 
$$
\delta^{r+1}(u)= \sum _{p=1}^{s}(-1)^{\mid \{q : i_{q}< j_{p}\}\mid }
\psi _{j_p}(u)\\
=\sum _{p=1}^{s}(-1)^{j_{p}-p}\psi _{j_p}(u).
$$

\par
For $F \subset [n]$, we define $x^{F}
:=\prod_{i \in F}  x_{i}$.
We define a simplicial complex 
$\Delta(F)$ 
by
\begin{displaymath}
  \Delta(F)=
  \big\{
    \{ i_{1},i_{2}, \dots , i_{r}\} \subset [\mu] \; : \ \ 
    x^{F}\mid 
    \prod  _{j \in \{ 1,2, \cdots \mu \}\setminus
    \{i_{1},i_{2}, \dots , i_{r}\}}m_{j} \big\}.
\end{displaymath}

\par
For $\alpha \in {\bf Z}^{n}$, there is a unique decomposition
$\alpha=\alpha _{+} - \alpha_{-}$ such that $\alpha_{+}, \alpha_{-} \in 
{\bf N}^{n}$ and 
$\mbox{\rm supp }\alpha_{+} \cap \mbox{\rm supp }\alpha_{-}
= \emptyset$.
Then we have $\mbox{\rm supp }\alpha_{-}=\{ i : \alpha_{i}<0 \}.$

\begin{lemma}[{\cite[Lemma 2.5]{St1}}]
For $\alpha \in {\bf Z}^{n}$
give an orientation for $
\Delta(\mbox{\rm supp }\alpha_{-})$
by $1 <2 < \cdots < \mu $.
Then we have the following isomorphism of complexes$:$ 
$C^{\bullet}_{\alpha} \cong \tilde{C}_{\bullet}
(\Delta(\mbox{\rm supp }\alpha_{-}))$
such that 
$C^{r}_{\alpha} \cong \tilde{C}_{\mu-r-1}
(\Delta(\mbox{\rm supp }\alpha_{-}))$.
\end{lemma}

\par
Then using the above lemma we have the following theorem:

\begin{theorem}
\label{claim:mu-1}
The following conditions are equivalent for a squarefree monomial ideal $I:$ 
\begin{enumerate}
\item[(1)] $\pd_S S/I = \mu (I)-1$. 
\item[(2)] The hypergraph $\calH:=\calH(I)$ satisfies $B(\calH) \neq V(\calH)$ and 
either one of the following conditions$:$ 
\begin{enumerate}
\item[(i)] The graph $(W(\calH), \calH_{W(\calH)}^{1})$ contains 
  a complete bipartite graph as a spanning subgraph. 
\item[(ii)] There exists $i \in B(\calH)$ such that 
  $\{i, j\} \in \calH$ for all $j \in W(\calH)$.
\end{enumerate}
\end{enumerate}
\end{theorem} 
\begin{proof}
We may assume $B(\calH) \ne V(\calH)$. 
By the above lemma we have the following isomorphisms:
\begin{eqnarray*}
H^{\mu -1}_{I}(S)_{\alpha}  
& = & 
H^{\mu -1}(C^{\bullet})_{\alpha} \\ 
& \cong & 
\tilde {H}_{0}(\Delta(\mbox{\rm supp }\alpha_{-}); K)\\
& \cong &
\tilde {H}_{0}(\Delta(\mbox{\rm supp }\alpha_{-})^{(1)}; K ),
\end{eqnarray*}
where $\Delta(\mbox{\rm supp }\alpha_{-})^{(1)}
:=\{ F \in \Delta(\mbox{\rm supp }\alpha_{-}): \  \dim F \le 1 \}$ is the 1-skeleton of
$\Delta(\mbox{\rm supp }\alpha_{-})$. 
Hence $H^{\mu -1}_{I}(S)_{\alpha} = 0$ if and only if 
$\Delta(\mbox{\rm supp }\alpha_{-})^{(1)}$ is connected.

\par
We claim that $\Delta(\mbox{\rm supp }\alpha_{-})^{(1)}$ is connected
for all $\alpha \in {\bf Z}^{n}$ if and only if 
the graph $(U, {U \choose 2} \setminus \calH^{1}_{U})$ 
is connected for all $W(\calH) \subset U \subset V(\calH)$,
where ${U \choose 2} :=\{ \{i,j \} \subset U: \ i \ne j \}$. 
Let $U$ be the vertex set of  
$\Delta(\mbox{\rm supp }\alpha_{-})^{(1)}$ for $\alpha \in {\bf Z}^{n}$.
Then $U \supset W(\calH)$ and   ${U \choose 2} \setminus \calH^{1}_{U} \subset  
\Delta(\mbox{\rm supp }\alpha_{-})^{(1)}$.
Hence if  $(U, {U \choose 2} \setminus \calH^{1}_{U})$  is connected,
then so is $\Delta(\mbox{\rm supp }\alpha_{-})^{(1)}$. 
On the other hand,
fix $U$ such that $W(\calH) \subset U \subset V(\calH)$.
Put $U=W(\calH) \cup B'$, where $B' \subset B(\calH)$.
By a suitable change of variables, we may assume that $B'=\{1,2,\dots ,p \}$.
For $1 \le j \le p$, set 
$\{x_{i}: \, x_{i}| m_{j},  x_{i} \nmid m_{\ell}
\mbox{ for } 1 \le \ell \le \mu \mbox{ with }  \ell \ne j\}
=\{x_{i_{j1}}, \dots , x_{i_{js_{j}}}\}$.
Take $\alpha \in {\bf Z}^{n}$ such that 
$\mbox{\rm supp }\alpha_{-}=[n]\setminus 
\{i_{11}, \dots , i_{1s_{1}}, i_{21}, \dots , i_{2s_{2}}, 
\dots , i_{p1}, \dots , i_{ps_{p}} \}$.
Then we have 
${U \choose 2} \setminus \calH^{1}_{U}$ is equal to the set of all 2-faces in
$\Delta(\mbox{\rm supp }\alpha_{-})$. Hence we have the claim.

\par
The graph $(U, {U \choose 2} \setminus \calH^{1}_{U})$ is connected for all $U$
such that $W(\calH) \subset U \subset V(\calH)$
if and only if the following conditions (I) and (II) are satisfied: 
\begin{enumerate}
\item[(I)] The graph 
  $(W(\calH), {W(\calH) \choose 2} \setminus \calH^{1}_{W(\calH) })$ 
  is connected.
\item[(II)] For $i \in B(\calH)$ set $U_{i}=W(\calH) \cup \{ i \}$.
  The graph $(U_{i}, {U_{i} \choose 2} \setminus \calH^{1}_{U_{i}})$ 
  is connected for all $i \in B(\calH)$. 
\end{enumerate}
Hence the condition  $\pd_S S/I = \mu (I)-1$ holds if and only if
one of the following conditions (i)' or (ii)' is satisfied: 
\begin{enumerate}
\item[(i)'] The graph 
  $(W(\calH), {W(\calH) \choose 2} \setminus \calH^{1}_{W(\calH) })$ 
  is disconnected.
\item[(ii)'] The graph 
  $(U_{i}, {U_{i} \choose 2} \setminus \calH^{1}_{U_{i}})$ is disconnected 
  for some $i \in B(\calH)$. 
\end{enumerate}
The condition (i)'  is equivalent 
to the condition (i).
Under the assumption that the condition (i) (or equivalently (i)') does not hold
the conditions (ii) and (ii)' are equivalent.
Hence we are done.
\end{proof}

\begin{corollary}
\label{claim:char-free}
The condition $\pd_S S/I = \mu (I)-1$ is independent of the base field $K$
for a monomial ideal $I$.
\end{corollary}

\section{Generic set}
We cannot restrict the number of variables in general when we classify squarefree monomial ideals 
with a certain given property.
In that case, it is convenient to consider finitely generated squarefree monomial ideals 
in a polynomial ring with infinite variables. 

\par
Let $K$ be a field and $S_{\infty} = K[x_1, x_2, \ldots]$ a 
polynomial ring with countably infinite variables over $K$. 
Let $\mathcal{I}$ be the set of all finitely generated squarefree monomial 
ideals of $S_{\infty}$. 
For $I \in \mathcal{I}$, we denote by $X(I)$ the set of all variables 
which appear in one of the minimal monomial generators of $I$. 

\begin{definition}
  \label{defn:GenSet}
  Let $\mathcal{C}$ be some property on $\mathcal{I}$. 
  That is, the subset 
  \begin{displaymath}
    \mathcal{I}(\mathcal{C}) 
    := \{ I \in \mathcal{I}: I \mbox{ satisfies the property } \mathcal{C}\}
  \end{displaymath}
  is uniquely determined.
  A subset $\mathcal{A} \subset \mathcal{I}$ is called a 
  \textit{generic set} on $\mathcal{C}$ if the following two conditions 
  are satisfied: 
  \begin{enumerate}
  \item $\mathcal{A} \subset \mathcal{I} (\mathcal{C})$. 
  \item For any $J \in \mathcal{I} (\mathcal{C})$ with 
   $G(J) = \{ m_1, \ldots, m_{\mu} \}$, 
    there exist $I \in \mathcal{A}$ with $X(I) = \{ x_1, \ldots, x_n \}$ and 
    $G(I) = \{ m_1', \ldots, m_{\mu}' \}$ 
    where $m_{i}' = x_{t_{i1}} \cdots x_{t_{i j_i}}$, 
    and (possibly trivial) squarefree monomials $M_1, \ldots, M_n$ 
    on $X(J)$ that are pairwise coprime 
    such that 
    \begin{displaymath}
      M_{t_{i1}} \cdots M_{t_{i j_i}} = m_i, 
      \qquad i = 1, 2, \ldots, \mu. 
    \end{displaymath}
  \end{enumerate}

  \par
  We say that a generic set $\mathcal{A}$ on $\mathcal{C}$ is \textit{minimal} 
  if $\mathcal{A}$ is minimal among generic sets on $\mathcal{C}$ 
  with respect to inclusion, and we say $\mathcal{A}$ on $\mathcal{C}$ 
  is \textit{reduced} if 
  $$
  K[x_F : F \in \calH (I)]/I_{\calH (I)} \cong K[X(I)]/ (I \cap K[X(I)])
  $$ 
  for all $I \in \mathcal{A}$. 
\end{definition}

\par
A minimal generic set has the following property: 
\begin{proposition}
  \label{claim:GenSet}
  Let $\mathcal{A}$ be a minimal generic set on some property $\mathcal{C}$. 
  \begin{enumerate}
  \item For $J \in \mathcal{I} (\mathcal{C})$, 
    there exists $I \in \mathcal{A}$ such that 
    $\calH (J) \subset \calH (I)$ and $\mu (J) = \mu (I)$. 
  \item If $\calH (I') \subset \calH (I)$ and $\mu (I') = \mu (I)$ 
    for $I', I \in \mathcal{A}$, then $I' = I$. 
  \end{enumerate}
\end{proposition}
\begin{proof}
  (1) By definition, there exists an ideal $I \in \mathcal{A}$ 
  with the condition (2) of Definition \ref{defn:GenSet}. 
  In particular, $\mu (I) = \mu (J) =: \mu$. 
  We use the same notations as in Definition \ref{defn:GenSet}. 
  Take a face $F \in \calH (J)$ and let $y_{F}$ be a defining variable of $F$ 
  in $\calH (J)$. Take $i \in F$. Then $y_F$ divides $m_{i}$. 
  Since $m_{i}$ can be written as the product 
  $M_{t_{i 1}} \cdots M_{t_{i j_i}}$, we may assume that $M_{t_{i 1}}$ 
  is divisible by $y_F$. Then for $1 \leq \ell \leq \mu$, the variable $y_F$ 
  divides $m_{\ell}$ if and only if the variable $x_{t_{i1}}$ divides 
  $m_{\ell}'$. This means that $x_{t_{i1}}$ defines $F$ of $\calH (I)$. 

  \par
  (2) First note that in the proof of (1), we also proved that 
  $\calH (J) \subset \calH (I)$ holds when $I$ and $J$ have the connection 
  as in Definition \ref{defn:GenSet} (2). 
  Therefore from the minimality of $\mathcal{A}$, 
  it is enough to prove that if $\calH (J) \subset \calH(I)$ with 
  $\mu (J) = \mu (I) = \mu$ 
  for $I, J \in \mathcal{I}$, then $I$ and $J$ have that connection. 

  \par
  Let $G(J) = \{ m_1, \ldots, m_{\mu} \}$ and 
  $G(I) = \{ m_1', \ldots, m_{\mu}' \}$. 
  For $F \in \calH (J)$ (resp.\  $G \in \calH (I)$), we denote by $M_F$ 
  (resp.\  $M_G'$), the product of all defining variables 
  of $F$ in $\calH (J)$ (resp.\  $G$ in $\calH (I)$). 
  Then 
  \begin{displaymath}
    m_i = \prod_{\genfrac{}{}{0pt}{}{F \in \calH (J)}{i \in F}} M_F, 
    \qquad 
    m '_i = \prod_{\genfrac{}{}{0pt}{}{F \in \calH (J)}{i \in F}} M_F' 
          \prod_{\genfrac{}{}{0pt}{}{G \in \calH (I) \setminus \calH (J)}{
                  i \in G}} M_G'. 
  \end{displaymath}
  Since $F \in \calH (I)$, the product $M_F' \neq 1$. 
  Let $x_F$ be a variable which divides $M_F'$. Then the substitution 
  $M_F$ to $x_F$; $1$ to the variables dividing one of $M_F'/x_F$, $M_G'$ 
  yields the desired connection.      
\end{proof}
 
By virtue of Propositions \ref{claim:hgraph} and \ref{claim:GenSet}
we have the following proposition: 
\begin{proposition}
  \label{claim:Reduction}
  Let $\mathcal{A}$ be a generic set on some property $\mathcal{C}$. 
 Suppose for all $I, J \in \mathcal{I}(\mathcal{C})$
 we have 
 $$
 \pd_{K[X(I)]} K[X(I)]/(I \cap K[X(I)]) =\pd _{K[X(J)]} K[X(J)]/(J \cap K[X(J)]).
 $$ 
 We also assume that 
 $$
 \ara (I \cap K[X(I)]) = \pd_{K[X(I)]} K[X(I)]/(I \cap K[X(I)]) 
 $$ 
 for all $I \in \mathcal{A}$. 
 Then we have 
$$
 \ara (I \cap K[X(I)]) = \pd_{K[X(I)]} K[X(I)]/(I \cap K[X(I)]) 
 $$ 
 for all $I \in \mathcal{I}(\mathcal{C})$.
\end{proposition}

\par
Let $I$ be a squarefree monomial ideal of $S$. 
The next proposition guarantees that we may assume that 
$\calH (I)$ is connected when we consider the problem on the arithmetical rank.
\begin{proposition}
  \label{claim:connected}
  Let $I_1, I_2$ be squarefree monomial ideals of $S$. 
  Suppose that $X (I_1) \cap X (I_2) = \emptyset$. Then we have 
  \begin{displaymath}
    \ara (I_1 + I_2) \leq \ara I_1 + \ara I_2. 
  \end{displaymath}

  \par
  Moreover, if $\ara I_i = \pd_S S/{I_i}$ holds for $i = 1, 2$, 
  then $\ara (I_1 + I_2) = \pd_S S/{(I_1 + I_2)}$ also holds. 
\end{proposition}
\begin{proof}
  When $f_1, \ldots, f_{s_1} \in S$ generate $I_1$ up to radical and 
  $g_1, \ldots, g_{s_2} \in S$ generate $I_2$ up to radical, 
  it is clear that these $s_1 + s_2$ elements 
  generate $I_1 + I_2$ up to radical. 
  Hence, we have the desired inequality. 

  \par
  Let $F_{\bullet}, G_{\bullet}$ be minimal free resolutions 
  of $S/I_1, S/I_2$, respectively. 
  Then $F_{\bullet} \otimes G_{\bullet}$ provides a minimal free resolution 
  of $S/{(I_1 + I_2)}$, and its length is equal to 
  \begin{displaymath}
    \pd_S S/I_1 + \pd_S S/I_2 = \ara I_1 + \ara I_2 \geq \ara (I_1 + I_2)  
  \end{displaymath} 
  when $\ara I_i = \pd_S S/{I_i}$ holds for $i = 1, 2$. 
  Since the inequality $\ara (I_1 + I_2) \geq \pd_S S/(I_1 + I_2)$ 
  always holds by (\ref{eq:ara>pd}), 
  we have the desired equality. 
\end{proof}

\par
Let $I$ be a squarefree monomial ideal 
with $G(I) = \{ m_1, \ldots, m_{\mu} \}$. 
We say that $I$ is \textit{connected} 
if for any distinct indices $i, j \in [\mu]$, 
there exists a sequence of indices $i=i_0, i_1, \ldots, i_{\ell} = j$ 
such that $\gcd (m_{i_{k-1}}, m_{i_k}) \neq 1$ for all $k=1, \ldots,\ell$. 
In other words, 
the corresponding hypergraph $\calH (I)$ is connected. 

\begin{remark}
  \label{rmk:connected}
  Proposition \ref{claim:connected} guarantees that 
  we may assume $\indeg I \geq 2$. 
  If $\indeg I = 1$, then $I$ can be written as $I_1 + (x)$, 
  where $x \notin X (I_1)$. 
  Then $\pd_S S/I = \pd_S S/{I_1} + 1$ and $\ara I \leq \ara I_1 + 1$. 
  Therefore if $\ara I_1 = \pd_S S/{I_1}$ holds, then $\ara I = \pd_S S/I$ 
  also holds. 
\end{remark}

\section{An algorithm to find a minimal reduced generic set}
\label{sec:algo}
In this section we present an algorithm for computing a minimal reduced 
generic set on the property: 
\begin{displaymath}
    {\mathcal{C}}: \  \mu (I) = 5, \  \pd_{K[{X(I)}]} K[{X(I)}]/(I \cap K[{X(I)}])  = 3, 
                      \  \height I = 2, \  \text{$I$ is connected}. 
\end{displaymath}

\par
The algorithm is composed mainly by three steps:
\begin{enumerate}
\item [Step 1)]
 List up all the connected hypergraphs with  five vertices which satisfy  the separability condition; 
\item [Step 2)]
  Choose the hypergraphs  $\calH$ with $\height I_{\calH}= 2$ and with 
   $\pd K[X(I_{\calH})]/I_{\calH}=3$;
\item [Step 3)]
  Find a minimal reduced generic set on ${\mathcal{C}}$.
\end{enumerate}

\par
Since we use a computer we must note that the projective dimension does not 
depend on the characteristic of $K$ when $\mu (I) = 5$. 

\begin{lemma}
  \label{claim:charfree}
  Let $I$ be a squarefree monomial ideal of $S$ with $\mu (I) = 5$. 
  The projective dimension $\pd_S S/I$ does not depend on the 
  characteristic of $K$. 
\end{lemma}
\begin{proof}
We fix the base field $K$.
If $\pd_S S/I =1$, then $I$ is a principal ideal and the claim is clear.

\par
Suppose $\pd_S S/I =2$. Then the height of $I$ is either 1 or 2.
We assume that $\height I=1$. Then the total Betti numbers do not change 
if we remove the prime components of height one from $I$.
Hence we may assume that $\height I=2$. 
Since $\pd_S S/I =\height I$, $S/I$ is Cohen-Macaulay.
Since Cohen-Macaulayness of the height two monomial ideals 
does not depend on the base field,
we are done in this case.

\par
If $\pd_S S/I \ge 4$, then $\pd_S S/I$ does not depend on the 
characteristic of $K$ by Proposition \ref{claim:mu} 
and Corollary \ref{claim:char-free}. 

\par
Finally, suppose $\pd_S S/I =3$. By the above argument
it is not possible to have $\pd_S S/I \ne 3$ for another base field. 
Then we are done.
\end{proof}

\par
Now we explain the details of each step. 

\par
\bigskip

\par
{\bf Step 1):} 
To list up all the connected hypergraphs with  five vertices,
we consider a correspondence between a hypergraph and a bipartite graph.

\par
First we give necessary definitions.

\par
Let $G$ be a graph on the vertex set $V$ with the edge set $E(G)$. 
For a vertex $v \in V$, the set of neighbourhoods of $v$ is defined by 
$N(v) := \{ u \in V : \{ u, v \} \in E(G) \}$. 
For a subset $V' \subset V$, we write $G_{V'}$ the induced subgraph of $G$ 
on $V'$. 
  For a vertex set $V$ we consider a labeled partition $V=X \cup Y$.
  It means that $X$ and $Y$ are labeled.
  We call $X$ the {\it first subset }(or the {\it indeterminate-part  subset}) of $V$ 
  and $Y$ the {\it second subset} (or the {\it generator-part  subset}) of $V$.
  Note that a bipartite graph with a labeled partition $V=X \cup Y$
  can be regarded as a directed bipartite graph on the partition $V=X \cup Y$
  such that every edge $\{x, y\}$ has the direction from $x\in X$ to $y \in Y$.
\begin{definition}
  \label{defn:incidence_graph}
  Set $X = \{ x_1, \ldots, x_n \}$, $Y = \{ y_1, \ldots, y_{\mu} \}$.
  A bipartite graph $G$ with the labeled partition $V = X \cup Y$ is said to be a 
  \textit{corresponding graph} if it satisfies the following three conditions:
  \begin{enumerate}
  \item The graph $G$ is connected;
  \item $N(y_i)\not \subset N(y_j)$ for all $i\neq j$; 
  \item $N(x_i) \ne N(x_j)$ for any  $i \neq j$. 
  \end{enumerate} 
\end{definition}

\par
 Two corresponding graphs $G$ and $G'$ with labeled partitions 
 $V=X \cup Y$ and $V'=X' \cup Y'$, respectively, are isomorphic 
 if there are bijections $f: X \longrightarrow X'$ and $g: Y \longrightarrow Y'$
 such that $\{x, y\} \in E(G)$ for $x\in X$ and $y \in Y$ if and only if  $\{f(x), g(y)\} \in E(G')$.  

\par
 Now we describe the correspondence more concretely.

\begin{proposition}
\label{claim:incidence_graph}
  We set $X = \{ x_1, \ldots, x_n \}$, $Y = \{ y_1, \ldots, y_{\mu} \}$ 
  and $V = X \cup Y$. 
  Let $\mathcal{G}$ be the set of all isomorphism classes of corresponding graphs 
with the labeled partition
  $V= X \cup Y$ and 
  $\mathcal{S}$ be the set of all isomorphism classes of connected  hypergraphs 
  with  the vertex set $Y$  and with $n$ faces 
  which satisfy  the separability condition.
\par
  Then the map $\phi : \mathcal{G} \rightarrow \mathcal{S}$ defined by 
  sending $G \in \mathcal{G}$ to 
  $\calH =\{F_1, \dots , F_n\} \in \mathcal{S}$ where  $F_i=\{y_j : \{x_i, y_j\} \in E(G)\}$
  gives a one-to-one correspondence between these two sets. 
\end{proposition}

\begin{proof}
 Let $G$ be a corresponding graph on the labeled partition $V = X \cup Y$. 
 The connectivity of $G$ corresponds to that of the hypergraph $\phi (G)$.
 The condition (2) of Definition \ref{defn:incidence_graph} corresponds to
 the separability condition for the hypergraph $\phi (G)$.
 Moreover the condition (3) of Definition \ref{defn:incidence_graph} 
 guarantees that the hypergraph $\phi (G)$ has $n$ faces. 
 Hence the map  $\phi$ is well defined.
 A hypergraph $\calH=\{F_1, \dots , F_n\}$ on the vertex set 
 $Y=\{y_1, \ldots, y_{\mu}\}$ 
 with the separability condition is associated to the bipartite graph 
 on the labeled partition $X \cup Y$, 
 where $X=\{x_1, \dots, x_n\}$
 with $E(G)=\{\{x_i,y_j\} : y_j \in F_i \}$.
 This correspondence gives the inverse map of $\phi$.
\end{proof}

\par
By virtue of Proposition \ref{claim:incidence_graph}, 
it is enough to enumerate all the corresponding graphs with labeled partition
$V= X \cup Y$ such that $\mid Y \mid =5$ and $\mid X \mid \le 2^{5}=32$
to list up the corresponding hypergraphs.

\par
To perform this computation we used an existent software called {\it Nauty }(see \cite{Mk}), 
whose main purpose is to calculate a set of non-isomorphic graphs 
in an efficient way.
We also customized it by a routine in C language to obtain graphs satisfying the conditions in Definition \ref{defn:incidence_graph}.
See the site \cite{KRT_soft} for technical information and examples of computation. 
In this site we also provide the source under GPL license to generate the software applications.

\par
As a result, in the case $\mu=5$ the number of the corresponding  graphs, 
or equivalently the number of the 5-vertex hypergraphs with the separability condition, 
is around $1.8\cdot 10^7$.
We denote the set of the above graphs by ${\mathcal{P}}_1$. 

\par
\bigskip

\par
{\bf Step 2):} 
To each element $G \in {\mathcal{P}}_1$,
we associate the squarefree monomial ideal $I=(m_1, \dots , m_5)$ in $S=K[x_1, \dots , x_{32}]$,
where $m_i=\prod_{x_j \in N(y_i)}x_j$ for $i=1,\dots,5$.

\par
Then we can use {\it CoCoA} (\cite{Co})
to determine whether $\height I =2 $ and $\pd S/I=3$.
Choose all graphs in ${\mathcal{P}}_1$  which  correspond to an ideal satisfying the above properties.
Set this subset as ${\mathcal{P}}_2$.

\par
By this computation we know that 
the cardinality of ${\mathcal{P}}_2$ is around $2.3 \cdot 10^6$.

\par
\bigskip

\par
{\bf Step 3): }
 For the last step of our algorithm 
 we need to find a minimal generic set on the property ${\mathcal{C}}$, or, 
 equivalently,
 the corresponding set of bipartite graphs.
 For this purpose we define a partial order on bipartite graphs:
\begin{definition}
  \label{defn:pord_biparGraph}
  Let $G$ and  $G'$  be bipartite graphs on the labeled partitions $X \cup Y$ 
  and $X' \cup Y$, respectively. 
  We introduce a partial order $\preceq$ by setting $G' \preceq G$ 
  if and only if there is a subset $X''$ of $X$ such that 
  $G_{X'' \cup Y} \cong G'$.
\end{definition}

\par
We want to list up all the maximal corresponding graphs in ${\mathcal{P}}_2$
with respect to the partial order $\preceq$.
To reach this goal we partition the set ${\mathcal{P}}_2$ of the corresponding graphs
\begin{displaymath}
 {\mathcal{P}}_2 = \bigcup {\mathcal{G}}_i
\end{displaymath}
where $G \in {\mathcal{G}}_i$ if the first subset $X$  has the cardinality $i$. 
Let $\ell$ be the maximum $i$ 
such that ${\mathcal{G}}_i \neq \emptyset$. 
Obviously each $G \in {\mathcal{G}}_{\ell}$ is a maximal graph. 
Set
${\mathcal{M}}_{\ell}={\mathcal{G}}_{\ell}$. 

\par
As the next step we list up the maximal graphs in ${\mathcal{G}}_{\ell - 1}$.
For each graph $G \in {\mathcal{M}}_{\ell}$ 
and for each vertex $x_i \in X $ 
let $G_i' = G_{V(G) \setminus {\{ x_i \}}}$. 
We make the list $\{G_i' :    G \in {\mathcal{G}}_{\ell},\   x_i \in X \}$.
If a graph in ${\mathcal{G}}_{\ell - 1}$ is isomorphic to some $G_i'$ 
then we remove it. At the end of this process we have only the maximal graphs 
in ${\mathcal{G}}_{\ell - 1}$.
We denote the set of the maximal graphs in ${\mathcal{G}}_{\ell - 1}$ by ${\mathcal{M}}_{\ell - 1}$.

\par
As the next step we similarly list up the maximal graphs in ${\mathcal{G}}_{\ell - 2}$.
For each graph $G \in {\mathcal{G}}_{\ell -1}$ 
and for each vertex $x_i \in  X$ 
let $G_i' = G_{V(G) \setminus {\{ x_i\}}}$. 
We make the list 
$
\{G_i' :  G \in {\mathcal{G}}_{\ell -1},\  x_i \in X \}.
$
If a graph in ${\mathcal{G}}_{\ell - 2}$ is isomorphic to a graph in the list, 
then we remove it. 
At the end of this process we have only the maximal graphs 
in ${\mathcal{G}}_{\ell - 2}$.
We denote the set of the maximal graphs in ${\mathcal{G}}_{\ell - 2}$ by ${\mathcal{M}}_{\ell - 2}$.

\par
We repeat a similar procedure for each $i \leq \ell -3$
and define ${\mathcal{M}}_i$ as above. 

\par
Set
\begin{displaymath}
 {\mathcal{P}}_3 = \bigcup {\mathcal{M}}_i.
\end{displaymath}
Then ${\mathcal{P}}_3 $ is the set of maximal corresponding graphs.
A corresponding set of squarefree monomial ideals gives a minimal generic set on $\mathcal{C}$.

\par
As the result performing the algorithm we have 
$\mid {\mathcal{P}}_3 \mid =3$.
More concretely, we have the following theorem:
\begin{theorem}
\label{claim:result}
Let $\mathcal{C}$ be the property on $ \mathcal{I}$ such that
\begin{displaymath}
  \mathcal{I}(\mathcal{C})
  =\left\{ I \in \mathcal{I} \; : \; 
  \begin{aligned}
    &\mu (I) = 5, \  \pd_{K[{X(I)}]} K[{X(I)}]/(I \cap K[{X(I)}]) = 3, \\
    &\height I = 2, \  \text{$I$ is connected} 
  \end{aligned}
  \right\}.
\end{displaymath}
Then 
$\{ J_1, J_2, J_3 \}$ is a minimal reduced generic set on $\mathcal{C}$,
where we define 
$J_i=(m_1, m_2, m_3, m_4, m_5 )$, $i=1,2,3$, as follows$:$

For $J_1$ where $\mid  X(J_1) \mid = 23:$ 
\begin{displaymath}
  \begin{aligned}
    m_1 &= x_1 x_3 x_7 x_8 x_9 x_{10} x_{12} x_{14} x_{15} x_{18} x_{19} 
           x_{20} x_{21} x_{22}, \\ 
    m_2 &= x_2 x_4 x_5 x_8 x_9 x_{11} x_{13} x_{14} x_{17} x_{18} x_{19} 
           x_{20} x_{21} x_{23}, \\
    m_3 &= x_3 x_5 x_9 x_{11} x_{12} x_{15} x_{16} x_{17} x_{19} 
           x_{20} x_{22} x_{23}, \\
    m_4 &= x_4 x_6 x_{10} x_{11} x_{12} x_{13} x_{16} x_{18} x_{19} 
           x_{21} x_{22} x_{23}, \\
    m_5 &= x_6 x_7 x_{10} x_{13} x_{14} x_{15} x_{16} x_{17} 
           x_{20} x_{21} x_{22} x_{23}. 
  \end{aligned}
\end{displaymath}

\par
For $J_2$ where $\mid  X(J_2) \mid  = 24:$ 
\begin{displaymath}
  \begin{aligned}
    m_1 &= x_1 x_4 x_7 x_9 x_{10} x_{12} x_{15} x_{16} x_{18} x_{19} 
           x_{20} x_{22} x_{23} x_{24}, \\
    m_2 &= x_2 x_5 x_8 x_9 x_{11} x_{14} x_{15} x_{17} x_{18} x_{19} 
           x_{21} x_{22} x_{23} x_{24}, \\
    m_3 &= x_3 x_6 x_7 x_8 x_{11} x_{12} x_{13} x_{16} x_{17} x_{19} 
           x_{20} x_{21} x_{23} x_{24}, \\
    m_4 &= x_4 x_5 x_{10} x_{11} x_{12} x_{13} x_{14} x_{15} 
           x_{20} x_{21} x_{22} x_{23},\\
    m_5 &= x_6 x_{10} x_{13} x_{14} x_{16} x_{17} x_{18} 
           x_{20} x_{21} x_{22} x_{24}. 
  \end{aligned}
\end{displaymath}

\par
For $J_3 $ where $\mid  X(J_3) \mid = 25:$ 
\begin{displaymath}
  \begin{aligned}
    m_1 &= x_1 x_5 x_6 x_7 x_{11} x_{12} x_{15} x_{18} x_{19} 
           x_{20} x_{21} x_{22} x_{23} x_{24}, \\
    m_2 &= x_2 x_5 x_8 x_9 x_{11} x_{13} x_{15} x_{16} x_{17} x_{18} 
           x_{21} x_{22} x_{24} x_{25}, \\ 
    m_3 &= x_3 x_6 x_8 x_{10} x_{12} x_{14} x_{15} x_{16} x_{17} 
           x_{20} x_{21} x_{23} x_{24} x_{25}, \\
    m_4 &= x_4 x_7 x_9 x_{10} x_{13} x_{14} x_{17} x_{18} x_{19} 
           x_{20} x_{22} x_{23} x_{24} x_{25}, \\
    m_5 &= x_{11} x_{12} x_{13} x_{14} x_{16} x_{19} 
           x_{21} x_{22} x_{23} x_{25}. 
  \end{aligned}
\end{displaymath}
\end{theorem}
  
\begin{remark}
We have performed a similar algorithm for 
a minimal generic set on the property: 
\begin{displaymath}
  {\mathcal{C}}: \  \mu (I) = 5, \  \pd_{K[{X(I)}]} K[{X(I)}]/(I \cap K[{X(I)}])= 3, 
                    \  \height I = 3, \  \text{$I$ is connected}. 
\end{displaymath}
As a result, we know that a minimal generic set consists of 9 ideals, which is
coincident with the list chosen from  all corresponding  hypergraphs in \cite{KTYdev2},
where a combinatorial argument is used without a computer.
\end{remark}

\section{The case $\mu (I) \leq 5$}
In this section, we compute the arithmetical rank of a squarefree monomial 
ideal $I$ with $\mu (I) \leq 5$. 
\begin{theorem}
  \label{claim:5gen}
  Let $I$ be a squarefree monomial ideal of $S=K[x_1, \ldots, x_n]$. 
  Suppose that $\mu (I) \leq 5$. 
  Then 
  \begin{displaymath}
    \ara I = \pd_S S/I. 
  \end{displaymath}
\end{theorem}
\begin{proof}
  If $\height I = 1$, then $I$ is of the form 
  \begin{displaymath}
    I = (x_1, \ldots, x_t)I'S, 
  \end{displaymath}
  where $I'$ is a squarefree monomial ideal of $S' = K[x_{t+1}, \ldots, x_n]$ 
  with $\height I' \geq 2$. 
  Then the equality $\ara I = \pd_S S/I$ holds 
  if the same equality holds for $I' \subset S'$. 
  Hence we may assume that $\height I \geq 2$. 
  By \cite{KTYdev1, KTYdev2}, 
  we have already known $\ara I = \pd_S S/I$ 
  when $\mu (I) - \height I \leq 2$ or $\mu (I) - \pd_S S/I \leq 1$. 
  The remain cases are that 
  $\mu (I) = 5$, $\height I = 2$, and $\pd_S S/I = 2, 3$.  
  If $\pd_S S/I = 2$, then $\pd_S S/I = \height I = 2$, 
  and thus $S/I$ is Cohen--Macaulay. 
  The equality $\ara I = \pd_S S/I$ for such an ideal was proved 
  in \cite{Kim_h2CM}. 
  Therefore we only need to consider the case that 
  $\mu (I) = 5$, $\height I = 2$, and $\pd_S S/I = 3$. 
  By Proposition \ref{claim:Reduction} and Theorem \ref{claim:result}, 
  it is enough to prove 
  $\ara I = \pd_S S/I$ for the following three cases: 
  \begin{enumerate}
  \item $I=J_1 \cap K[X]$ with $\mid  X \mid  = 23$;
  \item $I=J_2 \cap K[X]$ with $\mid X \mid = 24$;
  \item $I=J_3 \cap K[X]$ with $\mid  X \mid = 25$. 
  \end{enumerate}

  \par
  Since we have $\ara I \geq \pd_S S/I = 3$ by (\ref{eq:ara>pd}), 
  it is sufficient to show $\ara I \leq 3$. 

  \par
  \textbf{Case (1):} 
  We prove that the following $3$ elements $g_1, g_2, g_3$ 
  generate $I$ up to radical: 
  \begin{displaymath}
    g_1 = x_{16} f_1 f_2, 
    \quad g_2 = m_1 f_1 + m_4, 
    \quad g_3 = m_2 f_2 + m_5, 
  \end{displaymath}
  where 
  \begin{displaymath}
      f_1 = \gcd (m_1, m_3) + \gcd (m_4, m_5) m_2, \quad
      f_2 = \gcd (m_2, m_3) + \gcd (m_4, m_5) m_1. 
  \end{displaymath}
  Set $J = (g_1, g_2, g_3)$. 
  Note that 
  $\sqrt{x_{16} \gcd (m_1, m_3) \gcd (m_2, m_3)} = m_3$, 
  where $\sqrt{m}$ denotes a product of all variables those are appear 
  in a monomial $m$. 
  Then it is easy to see that $\sqrt{J} \subset I$. 
  We prove the opposite inclusion. 

  \par
  Since 
  $x_{16} m_4 m_5 = x_{16} (g_2 - f_1 m_1) (g_3 - f_2 m_2) \in J$, 
  we have $x_7 x_{14} x_{15} x_{17} x_{20} m_4 \in \sqrt{J}$ and 
  $x_4 x_{11} x_{12} x_{18} x_{19} m_5 \in \sqrt{J}$. 
  Then $x_{17} m_4^2 = x_{17} m_4 (g_2 - f_1 m_1) \in \sqrt{J}$ 
  because $x_7 x_{14} x_{15} x_{20}$ divides $m_1$. 
  Therefore we have $x_{17} m_4 \in \sqrt{J}$, and thus 
  $x_{17} f_1 m_1 = x_{17} (g_2 - m_4) \in \sqrt{J}$. 
  Since $x_{17}$ divides $\gcd (m_2, m_3)$, 
  the belonging $\gcd (m_2, m_3) g_1 \in J$ implies 
  $x_{16} \gcd (m_2, m_3) f_1 \in \sqrt{J}$. 
  Then $x_{16} \gcd (m_2, m_3) \gcd (m_1, m_3) \in \sqrt{J}$ and 
  $x_{16} \gcd (m_2, m_3) \gcd (m_4, m_5) m_2 \in \sqrt{J}$ 
  because $m_4$ divides $\gcd (m_1, m_3) \gcd (m_4, m_5) m_2$. 
  Hence $m_3, \gcd (m_4, m_5) m_2 \in \sqrt{J}$. 
  Similarly, we have $\gcd (m_4, m_5) m_1 \in \sqrt{J}$. 
  Then $g_2 \in J$ (resp.\  $g_3 \in J$) implies that $m_1, m_4 \in \sqrt{J}$ 
  (resp.\  $m_2, m_5 \in \sqrt{J}$), as desired. 

  \par
  \textbf{Case (2):}
  We prove that the following $3$ elements $g_1, g_2, g_3$ 
  generate $I$ up to radical: 
  \begin{displaymath}
    g_1 = x_{13} \left( \frac{m_4}{x_{13}} + m_3 \right), 
    \quad g_2 = m_1 m_2 \left( \frac{m_4}{x_{13}} + m_3 \right) + m_5, 
    \quad g_3 = m_1 + m_2 + m_3. 
  \end{displaymath}
  Set $J = (g_1, g_2, g_3)$. 
  Then it is clear that $\sqrt{J} \subset I$. 
  We prove the opposite inclusion. 

  \par
  Since $x_{13} m_5 = x_{13} (g_2 - m_1 m_2 ((m_4/x_{13}) + m_3)) 
  = x_{13} g_2 - m_1 m_2 g_1 \in J$, we have $m_5 \in \sqrt{J}$. 
  Then $x_{13} x_{18} m_3^2 = x_{18} m_3 (g_1 - m_4) \in \sqrt{J}$ 
  because $x_{18} m_3 m_4$ is divisible by $m_5 \in \sqrt{J}$. 
  This implies $ x_{18} m_3 \in \sqrt{J}$. 
  Since $x_{18}$ is a common factor of $m_1$ and $m_2$, 
  we have $m_3^2 = m_3 (g_3 - m_1 - m_2) \in \sqrt{J}$. 
  That is $m_3 \in \sqrt{J}$, and thus $m_4 = g_1 - x_{13} m_3, 
  m_1 + m_2 = g_3 - m_3, 
  m_1 m_2 (m_4/x_{13}) = g_2 - m_1 m_2 m_3 - m_5 \in \sqrt{J}$. 
  Since $m_1 m_2$ is divisible by $m_4/x_{13}$, 
  we have $m_1 m_2 \in \sqrt{J}$. 
  Therefore $m_1, m_2 \in \sqrt{J}$, as desired. 

  \par
  \textbf{Case (3):} 
  Since the length of the Lyubeznik resolution (\cite{Ly88}) of 
  $I = (m_5, m_1, m_2, m_3, m_4)$ with respect to this order of generators 
  is equal to $3$, the inequality $\ara I \leq 3$ follows 
  from \cite{Kim_lyures}. 
\end{proof}

\section{The case $\arithdeg I \leq 4$}
Let $I$ be a squarefree monomial ideal of $S$ with $\arithdeg I \leq 4$. 
In this section, we prove that $\ara I = \pd_S S/I$ holds 
for such an ideal $I$. 

\begin{theorem}
  \label{claim:Adual_4gen}
  Let $S$ be a polynomial ring over a field $K$ 
  and $I$ a squarefree monomial ideal of $S$ with $\arithdeg I \leq 4$. 
  Then 
  \begin{displaymath}
    \ara I = \pd_S S/I. 
  \end{displaymath}
\end{theorem}

\par
To prove the theorem, first, we reduce to simple cases. 
Then we find $\pd_S S/I$ elements which generate $I$ up to radical. 
On this step, we use the result due to Schmitt and Vogel \cite{SchmVo}. 
\begin{lemma}[{Schmitt and Vogel \cite[Lemma, p.\  249]{SchmVo}}]
  \label{claim:SchmVo}
  Let $R$ be a commutative ring and $\mathcal{P}$ a finite subset of $R$. 
  Suppose that subsets  
  ${\mathcal{P}}_0, {\mathcal{P}}_1, \ldots, {\mathcal{P}}_r$ 
  of $\mathcal{P}$ satisfy the following three conditions$:$ 
  \begin{enumerate}
  \item[{(SV1)}] $\bigcup_{\ell = 0}^r {\mathcal{P}}_{\ell} = \mathcal{P};$ 
  \item[{(SV2)}] $|{\mathcal{P}}_0| = 1;$ 
  \item[{(SV3)}] for all $\ell >0$ and for all elements 
    $a, a'' \in {\mathcal{P}}_{\ell}$ with $a \neq a''$, there exist an integer 
    ${\ell}' < \ell$ and an element $a' \in {\mathcal{P}}_{{\ell}'}$ 
    such that $a a'' \in (a')$. 
  \end{enumerate}
  Set $I = ({\mathcal{P}})$ and 
  \begin{displaymath}
    g_{\ell} = \sum_{a \in {\mathcal{P}}_{\ell}} a, 
    \quad \ell = 0, 1, \ldots, r. 
  \end{displaymath}
  Then $g_0, g_1, \ldots, g_r$ generate $I$ up to radical. 
\end{lemma}

\par
First, we reduce to the case where $\dim \calH (I^{\ast}) \leq 1$. 
To do this, we need the following lemma. 
\begin{lemma}
  \label{claim:red_to_1dim}
  Let $I = I_{\Delta}$ be a squarefree monomial ideal 
  of $S=K[x_1, \ldots, x_n]$ 
  with $\arithdeg I \leq 4$ and $\indeg I \geq 2$. 
  Then there exists a squarefree monomial ideal $I' = I_{{\Delta}'}$ 
  such that $\dim \calH (I'^{\ast}) \leq 1$ and that
  $\Delta$ is obtained from  ${\Delta}'$ adding cones recursively.  
\end{lemma}
\begin{proof}
  Since $\mu (I^{\ast}) \leq 4$ and $\height I^{\ast} \geq 2$, 
  we have $\dim \calH (I^{\ast}) \leq 2$ by \cite[Proposition 3.4]{KTYdev1}. 
  We prove the lemma by induction on the number $c$ of $2$-faces of 
  $\calH (I^{\ast})$. 

  \par
  If $c=0$, then $\dim \calH (I^{\ast}) \leq 1$ 
  and we may take ${\Delta}' = \Delta$. 

  \par
  Assume $c \geq 1$. Then $\arithdeg I = 4$ since $\indeg I \geq 2$, 
  and there exists a $2$-face $F \in \calH (I^{\ast})$. 
  We may assume $F = \{ 1, 2, 3 \}$ 
  and let $x_1, \ldots, x_t$ ($t \geq 1$) be all defining variables of $F$. 
  Let 
  \begin{displaymath}
    I = I_{\Delta} = P_1 \cap P_2 \cap P_3 \cap P_4
  \end{displaymath}
  be the minimal prime decomposition of $I$. 
  Then we can write $P_i = P_i' + (x_1)$ 
  for $i=1, 2, 3$, where $x_1 \notin P_i'$. 
  Note that $x_1 \notin P_4$. 
  Set 
  \begin{equation}
    \label{eq:prime_decomp}
      I' = P_1' \cap P_2' \cap P_3' \cap P_4 
      \quad (\subset S' = K[x_2, \ldots, x_n]), 
  \end{equation}
  and let ${\Delta}'$ be the simplicial complex associated with $I'$. 
  Then $\Delta = {\Delta}' \cup \cone_{x_1} (G_4)$, 
  where $G_4$ is the facet of $\Delta'$ corresponding to $P_4$. 
  If $t>1$, then $\calH ({I'}^{\ast}) = \calH (I^{\ast})$ 
  and (\ref{eq:prime_decomp}) is still the minimal prime decomposition 
  of $I'$. 
  Hence we can reduce the case $t=1$ by applying the above 
  operation $t-1$ times. 
  If (\ref{eq:prime_decomp}) is not a minimal prime decomposition of $I'$, 
  then $I' = P_1' \cap P_2' \cap P_3'$ 
  and $\dim \calH (I'^{\ast}) \leq 1$. 
  Otherwise $\calH (I'^{\ast}) = \calH (I^{\ast}) \setminus \{ F \}$ 
  and the number of $2$-feces of $\calH (I'^{\ast})$ is smaller than 
  that of $\calH (I^{\ast})$, as required. 
\end{proof}

\par
Barile and Terai \cite[Theorem 2, p.\  3694]{BariTera08} 
(see also \cite[Section 5]{Kim_h2CM}) proved that 
if $\ara I' = \pd_{S'} {S'}/{I'}$ holds, then $\ara I = \pd_{S} {S}/{I}$ 
also holds, 
where the notations are the same as in Lemma \ref{claim:red_to_1dim}. 
Therefore we only need to consider the case where 
$\dim \calH (I^{\ast}) \leq 1$. 

\par
Next, we reduce to the case where $\calH (I^{\ast})$ is connected. 
\begin{lemma}
  \label{claim:inter_section}
  Let $I_1, I_2$ be squarefree monomial ideals of $S$. 
  Suppose that $X (I_1) \cap X (I_2) = \emptyset$. 
  Then 
  \begin{displaymath}
    \ara (I_1 \cap I_2) \leq \ara I_1 + \ara I_2 - 1. 
  \end{displaymath}

  \par
  Moreover, if $\ara I_i = \pd_S S/{I_i}$ holds for $i = 1, 2$, 
  then $\ara (I_1 \cap I_2) = \pd_S S/{(I_1 \cap I_2)}$ also holds. 
\end{lemma}
\begin{proof}
  Let $\ara I_i = s_i + 1$ for $i=1, 2$. 
  Then there exist $s_1 + 1$ elements $f_0, f_1, \ldots, f_{s_1} \in S$ 
  (resp.\  $s_2 + 1$ elements $g_0, g_1, \ldots, g_{s_2} \in S$) 
  such that those generate $I_1$ (resp.\  $I_2$) up to radical. 

  \par
  Set 
  \begin{displaymath}
    h_{\ell} = \sum_{j=0}^{\ell} 
      f_{\ell - j} g_j, \quad \ell = 0, 1, \ldots, s_1 + s_2 
  \end{displaymath}
  and $J = (h_0, h_1, \ldots, h_{s_1 + s_2})$, 
  where $f_i = 0$ (resp.\ $g_j = 0$) if $i > s_1$ (resp.\  $j > s_2$).  
  If $\sqrt{J} = I_1 \cap I_2$, then 
  \begin{displaymath}
    \ara (I_1 \cap I_2) \leq s_1 + s_2 + 1 = \ara I_1 + \ara I_2 - 1. 
  \end{displaymath}

  \par
  We prove $\sqrt{J} = I_1 \cap I_2$. 
  Since $h_{\ell} \in I_1 \cap I_2$, we have $\sqrt{J} \subset I_1 \cap I_2$. 
  We prove the opposite inclusion. 
  Note that $I_1 \cap I_2 = I_1 I_2$ 
  since $X (I_1) \cap X (I_2) = \emptyset$. 
  By Lemma \ref{claim:SchmVo}, 
  we have $f_i g_j \in \sqrt{J}$ for all $0 \leq i \leq s_1$ and 
  for all $0 \leq j \leq s_2$. Then 
  \begin{displaymath}
    (f_0, f_1, \ldots, f_{s_1}) (g_0, g_1, \ldots, g_{s_2}) \subset \sqrt{J}. 
  \end{displaymath}
  Hence 
  \begin{displaymath}
    \sqrt{(f_0, f_1, \ldots, f_{s_1}) (g_0, g_1, \ldots, g_{s_2})} 
    \subset \sqrt{J}. 
  \end{displaymath}
  On the other hand,  
  \begin{displaymath}
    \begin{aligned}
      \sqrt{(f_0, f_1, \ldots, f_{s_1}) (g_0, g_1, \ldots, g_{s_2})}
      &\supset 
       \sqrt{(f_0, f_1, \ldots, f_{s_1})} \sqrt{(g_0, g_1, \ldots, g_{s_2})} \\
      &= I_1 I_2 = I_1 \cap I_2. 
    \end{aligned}
  \end{displaymath}
  Therefore we have $I_1 \cap I_2 \subset \sqrt{J}$, as desired. 

  \par
  To prove the second part of the lemma, we set $I = I_1 \cap I_2$. 
  Then $I^{\ast} = I_1^{\ast} + I_2^{\ast}$. 
  Note that $X (I_1) \cap X (I_2) = \emptyset$ implies 
  $X (I_1^{\ast}) \cap X (I_2^{\ast}) = \emptyset$. 
  Since we have already known the inequality $\ara I \geq \pd_S S/I$ 
  by (\ref{eq:ara>pd}), 
  it is sufficient to prove the opposite inequality. 

  \par
  By assumption, we have 
  \begin{equation}
    \label{eq:inter_section1}
    \ara I = \ara (I_1 \cap I_2) \leq \ara I_1 + \ara I_2 - 1 
                      = \pd_S S/I_1 + \pd_S S/I_2 - 1. 
  \end{equation}
  On the other hand, 
  since $X (I_1^{\ast}) \cap X (I_2^{\ast}) = \emptyset$, we have 
  \begin{displaymath}
    {\beta}_{i,j} (S/{I^{\ast}}) 
    = \sum_{\ell = 0}^j \sum_{m=0}^i 
        {\beta}_{m, \ell} (S/{I_1^{\ast}})
        {\beta}_{i-m, j - \ell} (S/{I_2^{\ast}}). 
  \end{displaymath}
  Assume that the regularity of $S/{I_1^{\ast}}$ (resp.\  $S/{I_2^{\ast}}$) 
  is given by ${\beta}_{i_1, j_1} (S/{I_1^{\ast}}) \neq 0$ 
  (resp.\  ${\beta}_{i_2, j_2} (S/{I_2^{\ast}}) \neq 0$). Then 
  \begin{displaymath}
    \reg S/{I_1^{\ast}} = j_1 - i_1, \qquad 
    \reg S/{I_2^{\ast}} = j_2 - i_2, 
  \end{displaymath}
  and ${\beta}_{i_1 + i_2, j_1 + j_2} (S/{I^{\ast}}) \geq 
  {\beta}_{i_1, j_1} (S/{I_1^{\ast}}) 
  {\beta}_{i_2, j_2} (S/{I_2^{\ast}}) > 0$. 
  Hence, 
  \begin{equation}
    \label{eq:inter_section2}
    \reg I^{\ast} 
    = \reg (S/{I^{\ast}}) + 1 \geq j_1 + j_2 - (i_1 + i_2) + 1. 
  \end{equation}
  Therefore 
  \begin{displaymath}
    \begin{alignedat}{3}
      \pd_S S/I = \reg I^{\ast} 
                &\geq \reg S/{I_1^{\ast}} + \reg S/{I_2^{\ast}} + 1 
                  &\quad &\text{(by (\ref{eq:inter_section2}))}\\
                &= \reg I_1^{\ast} + \reg I_2^{\ast} - 1 &\quad &\quad \\
                &= \pd_S S/{I_1} + \pd_S S/{I_2} -1 &\quad &\quad \\
                &\geq \ara I &\quad &\text{(by (\ref{eq:inter_section1}))}. 
    \end{alignedat}
  \end{displaymath} 
\end{proof}

\par
Let $I$ be a squarefree monomial ideal with $\arithdeg I \leq 4$. 
By Remark \ref{rmk:connected} we may assume $\indeg I \geq 2$. 
It was proved in \cite[Theorems 5.1 and 6.1]{KTYdev1} 
that $\ara I = \pd_S S/I$ holds 
when $\arithdeg I - \indeg I \leq 1$. 
Hence we only need to consider the case where 
$\arithdeg I =4$ and $\indeg I = 2$, 
equivalently, $\mu (I^{\ast}) = 4$ and $\height I^{\ast} = 2$. 
Then corresponding hypergraphs $\calH := \calH (I^{\ast})$ were classified in 
\cite[Section 3]{KTYdev2}. 
As a consequence of Lemmas \ref{claim:red_to_1dim} 
and \ref{claim:inter_section}, 
we may assume that $\calH$ is connected and $\dim \calH \leq 1$. 
Moreover since we know $\ara I = \pd_S S/I$ holds 
when $\arithdeg I = \reg I$, i.e., $\mu (I^{\ast}) = \pd_S S/I^{\ast}$ 
by \cite[Theorem 5.1]{KTYdev1}, 
we may assume that $W({\calH}) \neq \emptyset$ by Proposition \ref{claim:mu}. 
Thus $\calH$ coincides with one of the following hypergraphs: 

\par
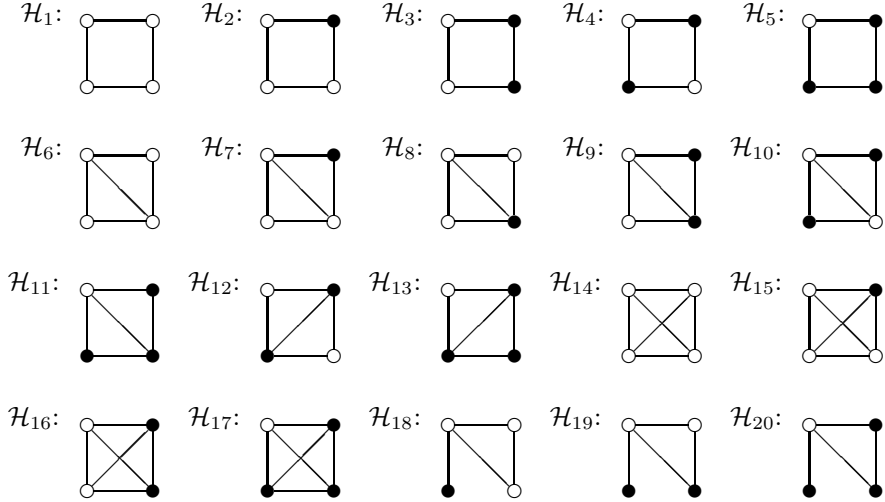
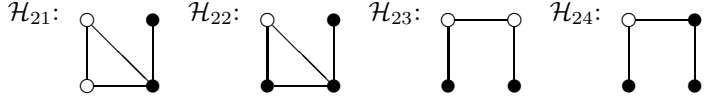
\begin{picture}(65,50)(-15,0)
    \put(-15,35){${\calH}_{1}$:}
    \put(10,10){\circle{5}}
    \put(35,10){\circle{5}}
    \put(10,35){\circle{5}}
    \put(35,35){\circle{5}}
    \put(12.5,10){\line(1,0){20}}
    \put(12.5,35){\line(1,0){20}}
    \put(10,12.5){\line(0,1){20}}
    \put(35,12.5){\line(0,1){20}}
\end{picture}
\begin{picture}(65,50)(45,0)
    \put(45,35){${\calH}_2$:}
    \put(70,10){\circle{5}}
    \put(95,10){\circle{5}}
    \put(70,35){\circle{5}}
    \put(95,35){\circle*{5}}
    \put(72.5,10){\line(1,0){20}}
    \put(72.5,35){\line(1,0){20}}
    \put(70,12.5){\line(0,1){20}}
    \put(95,12.5){\line(0,1){20}}
\end{picture}
\begin{picture}(65,50)(105,0)
    \put(105,35){${\calH}_3$:}
    \put(130,10){\circle{5}}
    \put(155,10){\circle*{5}}
    \put(130,35){\circle{5}}
    \put(155,35){\circle*{5}}
    \put(132.5,10){\line(1,0){20}}
    \put(132.5,35){\line(1,0){20}}
    \put(130,12.5){\line(0,1){20}}
    \put(155,12.5){\line(0,1){20}}
\end{picture}
\begin{picture}(65,50)(165,0)
    \put(165,35){${\calH}_4$:}
    \put(190,10){\circle*{5}}
    \put(215,10){\circle{5}}
    \put(190,35){\circle{5}}
    \put(215,35){\circle*{5}}
    \put(192.5,10){\line(1,0){20}}
    \put(192.5,35){\line(1,0){20}}
    \put(190,12.5){\line(0,1){20}}
    \put(215,12.5){\line(0,1){20}}
\end{picture}
\begin{picture}(65,50)(225,0)
    \put(225,35){${\calH}_5$:}
    \put(250,10){\circle*{5}}
    \put(275,10){\circle*{5}}
    \put(250,35){\circle{5}}
    \put(275,35){\circle*{5}}
    \put(252.5,10){\line(1,0){20}}
    \put(252.5,35){\line(1,0){20}}
    \put(250,12.5){\line(0,1){20}}
    \put(275,12.5){\line(0,1){20}}
\end{picture}

\par
\begin{picture}(65,50)(285,0)
    \put(285,35){${\calH}_6$:}
    \put(310,10){\circle{5}}
    \put(335,10){\circle{5}}
    \put(310,35){\circle{5}}
    \put(335,35){\circle{5}}
    \put(312.5,10){\line(1,0){20}}
    \put(312.5,35){\line(1,0){20}}
    \put(310,12.5){\line(0,1){20}}
    \put(335,12.5){\line(0,1){20}}
    \put(312,33){\line(1,-1){21}}
\end{picture}
\begin{picture}(65,50)(-15,0)
    \put(-15,35){${\calH}_7$:}
    \put(10,10){\circle{5}}
    \put(35,10){\circle{5}}
    \put(10,35){\circle{5}}
    \put(35,35){\circle*{5}}
    \put(12.5,10){\line(1,0){20}}
    \put(12.5,35){\line(1,0){20}}
    \put(10,12.5){\line(0,1){20}}
    \put(35,12.5){\line(0,1){20}}
    \put(12,33){\line(1,-1){21}}
\end{picture}
\begin{picture}(65,50)(45,0)
    \put(45,35){${\calH}_8$:}
    \put(70,10){\circle{5}}
    \put(95,10){\circle*{5}}
    \put(70,35){\circle{5}}
    \put(95,35){\circle{5}}
    \put(72.5,10){\line(1,0){20}}
    \put(72.5,35){\line(1,0){20}}
    \put(70,12.5){\line(0,1){20}}
    \put(95,12.5){\line(0,1){20}}
    \put(71.5,33.5){\line(1,-1){21.5}}
\end{picture}
\begin{picture}(65,50)(105,0)
    \put(105,35){${\calH}_9$:}
    \put(130,10){\circle{5}}
    \put(155,10){\circle*{5}}
    \put(130,35){\circle{5}}
    \put(155,35){\circle*{5}}
    \put(132.5,10){\line(1,0){20}}
    \put(132.5,35){\line(1,0){20}}
    \put(130,12.5){\line(0,1){20}}
    \put(155,12.5){\line(0,1){20}}
    \put(132,33){\line(1,-1){21.5}}
\end{picture}
\begin{picture}(65,50)(165,0)
    \put(160,35){${\calH}_{10}$:}
    \put(190,10){\circle*{5}}
    \put(215,10){\circle{5}}
    \put(190,35){\circle{5}}
    \put(215,35){\circle*{5}}
    \put(192.5,10){\line(1,0){20}}
    \put(192.5,35){\line(1,0){20}}
    \put(190,12.5){\line(0,1){20}}
    \put(215,12.5){\line(0,1){20}}
    \put(192,33){\line(1,-1){21}}
\end{picture}

\par
\begin{picture}(65,50)(225,0)
    \put(220,35){${\calH}_{11}$:}
    \put(250,10){\circle*{5}}
    \put(275,10){\circle*{5}}
    \put(250,35){\circle{5}}
    \put(275,35){\circle*{5}}
    \put(252.5,10){\line(1,0){20}}
    \put(252.5,35){\line(1,0){20}}
    \put(250,12.5){\line(0,1){20}}
    \put(275,12.5){\line(0,1){20}}
    \put(252,33){\line(1,-1){21.5}}
\end{picture}
\begin{picture}(65,50)(285,0)
    \put(280,35){${\calH}_{12}$:}
    \put(310,10){\circle*{5}}
    \put(335,10){\circle{5}}
    \put(310,35){\circle{5}}
    \put(335,35){\circle*{5}}
    \put(312.5,10){\line(1,0){20}}
    \put(312.5,35){\line(1,0){20}}
    \put(310,12.5){\line(0,1){20}}
    \put(335,12.5){\line(0,1){20}}
    \put(311.5,11.5){\line(1,1){22}}
\end{picture}
\begin{picture}(65,50)(-15,0)
    \put(-20,35){${\calH}_{13}$:}
    \put(10,10){\circle*{5}}
    \put(35,10){\circle*{5}}
    \put(10,35){\circle{5}}
    \put(35,35){\circle*{5}}
    \put(12.5,10){\line(1,0){20}}
    \put(12.5,35){\line(1,0){20}}
    \put(10,12.5){\line(0,1){20}}
    \put(35,12.5){\line(0,1){20}}
    \put(11.5,11.5){\line(1,1){22}}
\end{picture}
\begin{picture}(65,50)(45,0)
    \put(40,35){${\calH}_{14}$:}
    \put(70,10){\circle{5}}
    \put(95,10){\circle{5}}
    \put(70,35){\circle{5}}
    \put(95,35){\circle{5}}
    \put(72.5,10){\line(1,0){20}}
    \put(72.5,35){\line(1,0){20}}
    \put(70,12.5){\line(0,1){20}}
    \put(95,12.5){\line(0,1){20}}
    \put(72,12){\line(1,1){21}}
    \put(72,33){\line(1,-1){21}}
\end{picture}
\begin{picture}(65,50)(105,0)
    \put(100,35){${\calH}_{15}$:}
    \put(130,10){\circle{5}}
    \put(155,10){\circle{5}}
    \put(130,35){\circle{5}}
    \put(155,35){\circle*{5}}
    \put(132.5,10){\line(1,0){20}}
    \put(132.5,35){\line(1,0){20}}
    \put(130,12.5){\line(0,1){20}}
    \put(155,12.5){\line(0,1){20}}
    \put(132,12){\line(1,1){21}}
    \put(132,33){\line(1,-1){21}}
\end{picture}

\par
\begin{picture}(65,50)(165,0)
    \put(160,35){${\calH}_{16}$:}
    \put(190,10){\circle{5}}
    \put(215,10){\circle*{5}}
    \put(190,35){\circle{5}}
    \put(215,35){\circle*{5}}
    \put(192.5,10){\line(1,0){20}}
    \put(192.5,35){\line(1,0){20}}
    \put(190,12.5){\line(0,1){20}}
    \put(215,12.5){\line(0,1){20}}
    \put(192,12){\line(1,1){21}}
    \put(192,33){\line(1,-1){21.5}}
\end{picture}
\begin{picture}(65,50)(225,0)
    \put(220,35){${\calH}_{17}$:}
    \put(250,10){\circle*{5}}
    \put(275,10){\circle*{5}}
    \put(250,35){\circle{5}}
    \put(275,35){\circle*{5}}
    \put(252.5,10){\line(1,0){20}}
    \put(252.5,35){\line(1,0){20}}
    \put(250,12.5){\line(0,1){20}}
    \put(275,12.5){\line(0,1){20}}
    \put(251.5,11.5){\line(1,1){22}}
    \put(252,33){\line(1,-1){21.5}}
\end{picture}
\begin{picture}(65,50)(285,0)
    \put(280,35){${\calH}_{18}$:}
    \put(310,10){\circle*{5}}
    \put(335,10){\circle{5}}
    \put(310,35){\circle{5}}
    \put(335,35){\circle{5}}
    \put(312.5,35){\line(1,0){20}}
    \put(310,12.5){\line(0,1){20}}
    \put(335,12.5){\line(0,1){20}}
    \put(312,33){\line(1,-1){21}}
\end{picture}
\begin{picture}(65,50)(-15,0)
    \put(-20,35){${\calH}_{19}$:}
    \put(10,10){\circle*{5}}
    \put(35,10){\circle*{5}}
    \put(10,35){\circle{5}}
    \put(35,35){\circle{5}}
    \put(12.5,35){\line(1,0){20}}
    \put(10,12.5){\line(0,1){20}}
    \put(35,12.5){\line(0,1){20}}
    \put(12,33){\line(1,-1){21.5}}
\end{picture}
\begin{picture}(65,50)(45,0)
    \put(40,35){${\calH}_{20}$:}
    \put(70,10){\circle*{5}}
    \put(95,10){\circle*{5}}
    \put(70,35){\circle{5}}
    \put(95,35){\circle*{5}}
    \put(72.5,35){\line(1,0){20}}
    \put(70,12.5){\line(0,1){20}}
    \put(95,12.5){\line(0,1){20}}
    \put(72,33){\line(1,-1){21.5}}
\end{picture}

\par
\begin{picture}(65,50)(105,0)
    \put(100,35){${\calH}_{21}$:}
    \put(130,10){\circle{5}}
    \put(155,10){\circle*{5}}
    \put(130,35){\circle{5}}
    \put(155,35){\circle*{5}}
    \put(132.5,10){\line(1,0){20}}
    \put(130,12.5){\line(0,1){20}}
    \put(155,12.5){\line(0,1){20}}
    \put(132,33){\line(1,-1){21.5}}
\end{picture}
\begin{picture}(65,50)(165,0)
    \put(160,35){${\calH}_{22}$:}
    \put(190,10){\circle*{5}}
    \put(215,10){\circle*{5}}
    \put(190,35){\circle{5}}
    \put(215,35){\circle*{5}}
    \put(192.5,10){\line(1,0){20}}
    \put(190,12.5){\line(0,1){20}}
    \put(215,12.5){\line(0,1){20}}
    \put(192,33){\line(1,-1){21.5}}
\end{picture}
\begin{picture}(65,50)(225,0)
    \put(220,35){${\calH}_{23}$:}
    \put(250,10){\circle*{5}}
    \put(275,10){\circle*{5}}
    \put(250,35){\circle{5}}
    \put(275,35){\circle{5}}
    \put(252.5,35){\line(1,0){20}}
    \put(250,12.5){\line(0,1){20}}
    \put(275,12.5){\line(0,1){20}}
\end{picture}
\begin{picture}(65,50)(285,0)
    \put(280,35){${\calH}_{24}$:}
    \put(310,10){\circle*{5}}
    \put(335,10){\circle*{5}}
    \put(310,35){\circle{5}}
    \put(335,35){\circle*{5}}
    \put(312.5,35){\line(1,0){20}}
    \put(310,12.5){\line(0,1){20}}
    \put(335,12.5){\line(0,1){20}}
\end{picture}

\par
Note that $\calH$ is contained in the hypergraph $\calH_{17}$: 

\par
\noindent
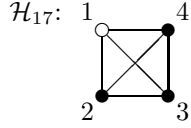
\begin{center}
  \begin{picture}(55,60)(-10,-5)
    \put(-25,39){${\calH}_{17}$:}
    \put(10,10){\circle*{5}}
    \put(2,1){$2$}
    \put(35,10){\circle*{5}}
    \put(38,1){$3$}
    \put(10,35){\circle{5}}
    \put(2,39){$1$}
    \put(35,35){\circle*{5}}
    \put(38,39){$4$}
    \put(12.5,10){\line(1,0){20}}
    \put(12.5,35){\line(1,0){20}}
    \put(10,12.5){\line(0,1){20}}
    \put(35,12.5){\line(0,1){20}}
    \put(11.5,11.5){\line(1,1){22}}
    \put(12,33){\line(1,-1){21.5}}
  \end{picture}
\end{center}

\par
\noindent
Throughout, we label the vertices of $\calH$ as above. 
Then $I$  is of the form 
$I = P_1 \cap P_2 \cap P_3 \cap P_4$ with 
\begin{displaymath}
  \begin{aligned}
    P_1 &= (x_{11}, \ldots, x_{1 i_1}, 
            x_{41}, \ldots, x_{4 i_4}, 
            x_{51}, \ldots, x_{5 i_5}), \\
    P_2 &= (x_{11}, \ldots, x_{1 i_1}, 
            x_{21}, \ldots, x_{2 i_2}, 
            x_{61}, \ldots, x_{6 i_6}, 
            y_{21}, \ldots, y_{2 j_2}), \\
    P_3 &= (x_{21}, \ldots, x_{2 i_2}, 
            x_{31}, \ldots, x_{3 i_3}, 
            x_{51}, \ldots, x_{5 i_5}, 
            y_{31}, \ldots, y_{3 j_3}), \\
    P_4 &= (x_{31}, \ldots, x_{3 i_3}, 
            x_{41}, \ldots, x_{4 i_4}, 
            x_{61}, \ldots, x_{6 i_6}, 
            y_{41}, \ldots, y_{4 j_4}), 
  \end{aligned}
\end{displaymath}
where $\{ x_{st} \}$, $\{ y_{u v} \}$ are all distinct variables of $S$ and 
$i_s \geq 0$, $j_u \geq 0$. 
Then we have 
\begin{displaymath}
  I^{\ast} = (X_1 X_4 X_5, X_1 X_2 X_6 Y_2, 
                 X_2 X_3 X_5 Y_3, X_3 X_4 X_6 Y_4), 
\end{displaymath}
where
\begin{displaymath}
  X_s = x_{s1} \cdots x_{s i_s}, \quad s = 1, 2, 3, 4; 
  \qquad
  Y_u = y_{u1} \cdots y_{u j_u}, \quad u = 2, 3, 4. 
\end{displaymath}
Here we set $X_s = 1$ (resp.\  $Y_u = 1$) when $i_s = 0$ (resp.\  $j_u=0$). 

\par
Let $N = i_1 + \cdots + i_6 + j_2 + j_3 + j_4$. 
Then one can easily construct a graded minimal free resolution of $I^{\ast} $ and 
compute $\reg I^{\ast}$. 
\begin{lemma}
  \label{claim:reg4gen}
  Let $I$ be a squarefree monomial ideal with $\indeg I^{\ast} \geq 2$. 
  Suppose that $\calH := \calH (I^{\ast})$ coincides with one of 
  $\calH_1, \ldots, \calH_{24}$. 
  \begin{enumerate}
  \item When $\calH$ coincides with one of 
    $\calH_{11}, \calH_{17}, \calH_{20}$, 
    \begin{displaymath}
      \pd_S S/I = \reg I^{\ast} 
      = \max \{ N - j_2 - 2, N - j_3 - 2, N - j_4 - 2 \}. 
    \end{displaymath}
  \item When $\calH = \calH_{22}$, 
    \begin{displaymath}
      \pd_S S/I = \reg I^{\ast} 
      = \max \{ N - j_2 - 2, N - j_3 - 2 \}. 
    \end{displaymath}
  \item When $\calH$ coincides with one of 
    $\calH_4, \calH_5, \calH_{12}, \calH_{13}, \calH_{24}$, 
    \begin{displaymath}
      \pd_S S/I = \reg I^{\ast} 
      = \max \{ N - j_2 - 2, N - j_4 - 2 \}. 
    \end{displaymath}
  \item Otherwise, $\pd_S S/I = \reg I^{\ast} = N-2$. 
  \end{enumerate}
\end{lemma}

\par
Now we find $\pd_S S/I$ elements of $I$ which generate $I$ up to radical. 
First, we consider the case of $\calH_{17}$. 
The following construction for $\calH = \calH_{17}$ is also valid 
for the other cases except for the two cases of $\calH_{1}$ and $\calH_{14}$. 

\par
Set 
\begin{displaymath}
  r_1 = N - j_4 - 3, \quad r_2 = N - j_3 -3, \quad r_3 = N - j_2 -3, 
  \quad r = \max \{ r_1, r_2, r_3 \}. 
\end{displaymath}
Then $\pd_S S/I = r+1$. 
We define sets 
${\mathcal{P}}_{\ell}^{(1)}, \ldots, {\mathcal{P}}_{\ell}^{(8)}$ 
as follows and set 
${\mathcal{P}}_{\ell} = {\mathcal{P}}_{\ell}^{(1)} \cup \cdots \cup
{\mathcal{P}}_{\ell}^{(8)}$ for $\ell = 0, 1, \ldots, r$: 
\begin{displaymath}
  \begin{aligned}
    {\mathcal{P}}_{\ell}^{(1)} 
    &= \left\{ x_{1 {\ell}_1} x_{3 {\ell}_3} \; : \; 
      {\ell}_1 + {\ell}_3 = \ell + 2; \  
      1 \leq {\ell}_s \leq i_s \  (s = 1, 3)
    \right\}, \\
    {\mathcal{P}}_{\ell}^{(2)} 
    &= \left\{ x_{1 {\ell}_1} w_{3 {\ell}_3} w_{4 {\ell}_4} \; : \; 
    \begin{aligned}
      &{\ell}_1 + {\ell}_3 + {\ell}_4 + i_3 = \ell + 3, \\
      &1 \leq {\ell}_1 \leq i_1; \  
      1 \leq {\ell}_3 \leq i_2 + i_5 + j_3; \  
      1 \leq {\ell}_4 \leq i_4 + i_6 + j_4
    \end{aligned}
    \right\}, \\
    {\mathcal{P}}_{\ell}^{(3)} 
    &= \left\{ x_{3 {\ell}_3} w_{1 {\ell}_1} w_{2 {\ell}_2} \; : \; 
    \begin{aligned}
      &{\ell}_3 + {\ell}_1 + {\ell}_2 + i_1 = \ell + 3, \\
      &1 \leq {\ell}_3 \leq i_3; \  
      1 \leq {\ell}_1 \leq i_4 + i_5; \ 
      1 \leq {\ell}_2 \leq i_2 + i_6 + j_2
    \end{aligned}
    \right\}, \\
    {\mathcal{P}}_{\ell}^{(4)} 
    &= \left\{ x_{2 {\ell}_2} x_{4 {\ell}_4} \; : \; 
      {\ell}_2 + {\ell}_4 + i_1 + i_3 = \ell + 2; \  
      1 \leq {\ell}_s \leq i_s \  (s = 2, 4)
    \right\}, \\
    {\mathcal{P}}_{\ell}^{(5)} 
    &= \left\{ x_{4 {\ell}_4} w_{2 {\ell}_2} w_{3 {\ell}_3} \; : \; 
    \begin{aligned}
      &{\ell}_4 + ({\ell}_2 - i_2) + ({\ell}_3 - i_2) + i_1 + i_2 + i_3 
        = \ell + 3, \\
      &1 \leq {\ell}_4 \leq i_4; \ 
      i_2 < {\ell}_2 \leq i_2 + i_6 + j_2; \ 
      i_2 < {\ell}_3 \leq i_2 + i_5 + j_3
    \end{aligned}
    \right\}, \\
    {\mathcal{P}}_{\ell}^{(6)} 
    &= \left\{ x_{2 {\ell}_2} x_{5 {\ell}_5} w_{4 {\ell}_4} \; : \; 
    \begin{aligned}
      &{\ell}_2 + {\ell}_5 + ({\ell}_4 - i_4) + i_1 + i_3 + i_4 = \ell + 3, \\
      &1 \leq {\ell}_s \leq i_s \  (s = 2, 5); \ 
      i_4 < {\ell}_4 \leq i_4 + i_6 + j_4
    \end{aligned}
    \right\}, \\
    {\mathcal{P}}_{\ell}^{(7)} 
    &= \left\{ x_{5 {\ell}_5} x_{6 {\ell}_6} \; : \; 
      {\ell}_5 + {\ell}_6 + i_1 + \cdots + i_4 = \ell + 2; \  
      1 \leq {\ell}_s \leq i_s \  (s = 5, 6)
    \right\}, \\
    {\mathcal{P}}_{\ell}^{(8)} 
    &= \left\{ x_{5 {\ell}_5} y_{2 {\ell}_2} y_{4 {\ell}_4} \; : \; 
    \begin{aligned}
      &{\ell}_5 + {\ell}_2 + {\ell}_4 + i_1 + \cdots + i_4 + i_6 = \ell + 3, \\
      &1 \leq {\ell}_5 \leq i_5; \ 
      1 \leq {\ell}_u \leq j_u \  (u = 2, 4) 
    \end{aligned}
    \right\}. 
  \end{aligned}
\end{displaymath}
Here, 
\begin{displaymath}
  \begin{aligned}
    w_{1 {\ell}_1} &= \left\{ 
    \begin{alignedat}{3}
      &x_{4 {\ell}_1}, &\quad & 1 \leq {\ell}_1 \leq i_4, \\
      &x_{5 {\ell}_1 - i_4}, &\quad & i_4 < {\ell}_1 \leq i_4 + i_5, 
    \end{alignedat}
    \right. \\
    w_{2 {\ell}_2} &= \left\{ 
    \begin{alignedat}{3}
      &x_{2 {\ell}_2}, &\quad & 1 \leq {\ell}_2 \leq i_2, \\
      &x_{6 {\ell}_2 - i_2}, &\quad & i_2 < {\ell}_2 \leq i_2 + i_6, \\
      &y_{2 {\ell}_2 - i_2 - i_6}, 
        &\quad & i_2 + i_6 < {\ell}_2 \leq i_2 + i_6 + j_2, 
    \end{alignedat}
    \right. \\
    w_{3 {\ell}_3} &= \left\{ 
    \begin{alignedat}{3}
      &x_{2 {\ell}_3}, &\quad & 1 \leq {\ell}_3 \leq i_2, \\
      &x_{5 {\ell}_3 - i_2}, &\quad & i_2 < {\ell}_3 \leq i_2 + i_5, \\
      &y_{3 {\ell}_3 - i_2 - i_5}, 
        &\quad & i_2 + i_5 < {\ell}_3 \leq i_2 + i_5 + j_3, 
    \end{alignedat}
    \right. \\
    w_{4 {\ell}_4} &= \left\{ 
    \begin{alignedat}{3}
      &x_{4 {\ell}_4}, &\quad & 1 \leq {\ell}_4 \leq i_4, \\
      &x_{6 {\ell}_4 - i_4}, &\quad & i_4 < {\ell}_4 \leq i_4 + i_6, \\
      &y_{4 {\ell}_4 - i_4 - i_6}, 
        &\quad & i_4 + i_6 < {\ell}_4 \leq i_4 + i_6 + j_4. 
    \end{alignedat}
    \right. 
  \end{aligned}
\end{displaymath}

\par
Note that for each $k$ $(k=1, 2, \ldots, 8)$, 
the range of $\ell$ with
${\mathcal{P}}_{\ell}^{(k)} \neq \emptyset$ is given by the 
following list: 
\begin{displaymath}
  \begin{aligned}
    {\mathcal{P}}_{\ell}^{(1)}: \quad &0 \leq \ell \leq i_1 + i_3 -2, \\
    {\mathcal{P}}_{\ell}^{(2)}: \quad &i_3 \leq \ell \leq N - j_2 - 3 = r_3, \\
    {\mathcal{P}}_{\ell}^{(3)}: \quad &i_1 \leq \ell \leq N - j_3 - j_4 - 3, \\
    {\mathcal{P}}_{\ell}^{(4)}: \quad 
       &i_1 + i_3 \leq \ell \leq i_1 + \cdots + i_4 -2, \\
    {\mathcal{P}}_{\ell}^{(5)}: \quad 
       &i_1 + i_2 + i_3 \leq \ell \leq N - j_4 - 3 = r_1, \\
    {\mathcal{P}}_{\ell}^{(6)}: \quad 
       &i_1 + i_3 + i_4 \leq \ell \leq N - j_2 - j_3 -3, \\
    {\mathcal{P}}_{\ell}^{(7)}: \quad 
       &i_1 + \cdots + i_4 \leq \ell \leq i_1 + \cdots + i_6 -2, \\
    {\mathcal{P}}_{\ell}^{(8)}: \quad 
       &i_1 + \cdots + i_4 + i_6 \leq \ell \leq N - j_3 - 3 = r_2. 
  \end{aligned}
\end{displaymath}

\par
Now, we verify that 
${\mathcal{P}}_{\ell}$, $\ell = 0, 1, \ldots, r$ 
satisfy the conditions (SV1), (SV2), and (SV3). 
In this case, the condition (SV1) means that 
$\bigcup_{\ell = 0}^{r} {\mathcal{P}}_{\ell}$ 
generates $I$, and it is satisfied. 
Since ${\mathcal{P}}_0 = {\mathcal{P}}_0^{(1)} = \{ x_{11} x_{31} \}$, 
the condition (SV2) is also satisfied. 
To check the condition (SV3), let $a, a''$ be two distinct elements in 
${\mathcal{P}}_{\ell}$.  
We denote the indices of $a$ (resp.\  $a''$) by ${\ell}_s$ 
(resp.\  ${\ell}_s''$). 
First suppose $a, a'' \in {\mathcal{P}}_{\ell}^{(k)}$. 
Then there exists $s$ such that ${\ell}_s \neq {\ell}_s''$ 
and we may assume ${\ell}_s < {\ell}_s''$. Replacing ${\ell}_s''$ 
to ${\ell}_s$ in $a''$, we obtain required elements 
$a' \in {\mathcal{P}}_{{\ell}'}^{(k)}$ with ${\ell}' < \ell$. 
For example, take two distinct elements 
$a = x_{1 {\ell}_1} x_{3 {\ell}_3}, 
a'' = x_{1 {\ell}_1''} x_{3 {\ell}_3''} \in {\mathcal{P}}_{\ell}^{(1)}$ 
with ${\ell}_1 < {\ell}_1''$. 
Then $a' = x_{1 {\ell}_1} x_{3 {\ell}_3''} \in {\mathcal{P}}_{{\ell}'}^{(1)}$, 
where ${\ell}' = {\ell}_1 + {\ell}_3'' - 2 < {\ell}_1'' + {\ell}_3'' - 2 
= \ell$. 
Next, we assume $a \in {\mathcal{P}}_{\ell}^{(k)}$ 
and $a'' \in {\mathcal{P}}_{\ell}^{(k'')}$ with $k < k''$. 
If $k=1$, $k'' = 2$, then $a = x_{1 {\ell}_1} x_{3 {\ell}_3}$, 
$a'' = x_{1 {\ell}_1''} w_{3 {\ell}_3''} w_{4 {\ell}_4''}$. 
We can take 
$a' = x_{1 {\ell}_1''} x_{3 {\ell}_3} \in {\mathcal{P}}_{{\ell}'}^{(1)}$, 
where 
\begin{displaymath}
  \begin{aligned}
    {\ell}' &= {\ell}_1'' + {\ell}_3 - 2 \\
            &= (\ell + 3 - ({\ell}_3'' + {\ell}_4'' + i_3)) + {\ell}_3 - 2 \\
            &= \ell + 1 - {\ell}_3'' - {\ell}_4'' - (i_3 - {\ell}_3) \\
            &< \ell. 
  \end{aligned}
\end{displaymath}
Similarly, we can check the existence of required ${\ell}'$ and 
$a' \in {\mathcal{P}}_{{\ell}'}^{(k')}$ for other pairs $(k, k'')$; 
we just mention the choices for $a'$ in Table \ref{tab:H1_list}. 

\begin{table}
  \caption{The choices for $a'$ in the case of ${\calH}_{17}$.}
  \label{tab:H1_list}
  \begin{center}
    \begin{tabular}{cccl} \hline
      $k$ & $k''$ & $k'$ & $a'$ \\ \hline
      $1$ & $2$   & $1$  & $x_{1 {\ell}_1''} x_{3 {\ell}_3}$ \\
      $1$ & $3$   & $1$  & $x_{1 {\ell}_1} x_{3 {\ell}_3''}$ \\
      $1$ & $4$--$8$ & \multicolumn{2}{c}{These cases do not occur.} \\
      $2$ & $3$   & $1$  & $x_{1 {\ell}_1} x_{3 {\ell}_3''}$ \\
      $2$ & $4$   & $2$  & 
        $x_{1 {\ell}_1} x_{2 {\ell}_2''} x_{4 {\ell}_4''} 
         = x_{1 {\ell}_1} w_{3 {\ell}_2''} w_{4 {\ell}_4''}$ \\
      $2$ & $5$   & $2$  & 
        $x_{1 {\ell}_1} w_{3 {\ell}_3''} x_{4 {\ell}_4''} 
         = x_{1 {\ell}_1} w_{3 {\ell}_3''} w_{4 {\ell}_4''}$ \\
      $2$ & $6$   & $2$  & 
        $x_{1 {\ell}_1} x_{2 {\ell}_2''} w_{4 {\ell}_4''} 
         = x_{1 {\ell}_1} w_{3 {\ell}_2''} w_{4 {\ell}_4''}$ \\
      $2$ & $7$   & $2$  & 
        $x_{1 {\ell}_1} x_{5 {\ell}_5''} x_{6 {\ell}_6''} 
         = x_{1 {\ell}_1} w_{3 i_2 + {\ell}_5''} w_{4 i_4 + {\ell}_6''}$ \\
      $2$ & $8$   & $2$  & 
        $x_{1 {\ell}_1} x_{5 {\ell}_5''} y_{4 {\ell}_4''} 
         = x_{1 {\ell}_1} w_{3 i_2 + {\ell}_5''} 
           w_{4 i_4 + i_6 + {\ell}_4''}$ \\
      $3$ & $4$   & $3$  & 
        $x_{3 {\ell}_3} x_{4 {\ell}_4''} x_{2 {\ell}_2''} 
         = x_{3 {\ell}_3} w_{1 {\ell}_4''} w_{2 {\ell}_2''}$ \\
      $3$ & $5$   & $3$  & 
        $x_{3 {\ell}_3} x_{4 {\ell}_4''} w_{2 {\ell}_2''} 
         = x_{3 {\ell}_3} w_{1 {\ell}_4''} w_{2 {\ell}_2''}$ \\
      $3$ & $6$   & $3$  & 
        $x_{3 {\ell}_3} x_{5 {\ell}_5''} x_{2 {\ell}_2''} 
         = x_{3 {\ell}_3} w_{1 i_4 + {\ell}_5''} w_{2 {\ell}_2''}$ \\
      $3$ & $7$   & $3$  & 
        $x_{3 {\ell}_3} x_{5 {\ell}_5''} x_{6 {\ell}_6''} 
         = x_{3 {\ell}_3} w_{1 i_4 + {\ell}_5''} w_{2 i_2+ {\ell}_6''}$ \\
      $3$ & $8$   & $3$  & 
        $x_{3 {\ell}_3} x_{5 {\ell}_5''} y_{2 {\ell}_2''} 
         = x_{3 {\ell}_3} w_{1 i_4 + {\ell}_5''} 
           w_{2 i_2 + i_6 + {\ell}_2''}$ \\
      $4$ & $5$   & $4$  & 
        $x_{2 {\ell}_2} x_{4 {\ell}_4''}$ \\
      $4$ & $6$   & $4$  & 
        $x_{2 {\ell}_2''} x_{4 {\ell}_4}$ \\
      $4$ & $7,8$ & \multicolumn{2}{c}{These cases do not occur.} \\
      $5$ & $6$   & $4$  & 
        $x_{2 {\ell}_2''} x_{4 {\ell}_4}$ \\
      $5$ & $7$   & $5$  & 
        $x_{4 {\ell}_4} x_{6 {\ell}_6''} x_{5 {\ell}_5''} 
         = x_{4 {\ell}_4} w_{2 i_2 + {\ell}_6''} w_{3 i_2+ {\ell}_5''}$ \\
      $5$ & $8$   & $5$  & 
        $x_{4 {\ell}_4} y_{2 {\ell}_2''} x_{5 {\ell}_5''} 
         = x_{4 {\ell}_4} w_{2 i_2 + i_6 + {\ell}_2''} 
           w_{3 i_2 + {\ell}_5''}$ \\
      $6$ & $7$   & $6$  & 
        $x_{2 {\ell}_2} x_{5 {\ell}_5''} x_{6 {\ell}_6''} 
         = x_{2 {\ell}_2} x_{5 {\ell}_5''} w_{4 i_4+ {\ell}_6''}$ \\
      $6$ & $8$   & $6$  & 
        $x_{2 {\ell}_2} x_{5 {\ell}_5''} y_{4 {\ell}_4''} 
         = x_{2 {\ell}_2} x_{5 {\ell}_5''} 
           w_{4 i_4 + i_6 + {\ell}_4''}$ \\
      $7$ & $8$   & $7$  & 
        $x_{5 {\ell}_5''} x_{6 {\ell}_6}$ \\ \hline
    \end{tabular}
  \end{center}
\end{table}

\par
\bigskip

\par
Next we consider the two exceptional cases: $\calH = \calH_1, \calH_{14}$: 

\par
\noindent
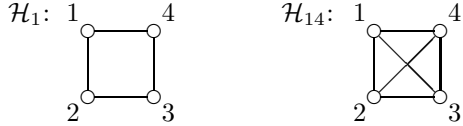
\begin{center}
  \begin{picture}(100,60)(-20,-5)
    \put(-20,39){${\calH}_1$:}
    \put(10,10){\circle{5}}
    \put(2,1){$2$}
    \put(35,10){\circle{5}}
    \put(38,1){$3$}
    \put(10,35){\circle{5}}
    \put(2,39){$1$}
    \put(35,35){\circle{5}}
    \put(38,39){$4$}
    \put(12.5,10){\line(1,0){20}}
    \put(12.5,35){\line(1,0){20}}
    \put(10,12.5){\line(0,1){20}}
    \put(35,12.5){\line(0,1){20}}
  \end{picture}
  \begin{picture}(90,60)(95,-5)
    \put(95,39){${\calH}_{14}$:}
    \put(130,10){\circle{5}}
    \put(122,1){$2$}
    \put(155,10){\circle{5}}
    \put(158,1){$3$}
    \put(130,35){\circle{5}}
    \put(122,39){$1$}
    \put(155,35){\circle{5}}
    \put(158,39){$4$}
    \put(132.5,10){\line(1,0){20}}
    \put(132.5,35){\line(1,0){20}}
    \put(130,12.5){\line(0,1){20}}
    \put(155,12.5){\line(0,1){20}}
    \put(132,12){\line(1,1){21}}
    \put(132,33){\line(1,-1){21}}
  \end{picture}
\end{center}

\par
The case of  ${\calH}_1$ is easy because in this case, 
\begin{displaymath}
  \begin{aligned}
    I &= (x_{11}, \ldots, x_{1 i_1}, x_{41}, \ldots, x_{4 i_4}) 
      \cap (x_{11}, \ldots, x_{1 i_1}, x_{21}, \ldots, x_{2 i_2}) \\
      &\cap (x_{21}, \ldots, x_{2 i_2}, x_{31}, \ldots, x_{3 i_3}) 
      \cap (x_{31}, \ldots, x_{3 i_3}, x_{41}, \ldots, x_{4 i_4}) \\
      &= (x_{1 {\ell}_1} x_{3 {\ell}_3} \; : \; 
          1 \leq {\ell}_s \leq i_s \  (s= 1, 3)) 
       + (x_{2 {\ell}_2} x_{4 {\ell}_4} \; : \; 
          1 \leq {\ell}_s \leq i_s \  (s= 2, 4)) \\ 
      &= (x_{11}, \ldots, x_{1 i_1}) 
      \cap (x_{31}, \ldots, x_{3 i_3}) 
       + (x_{21}, \ldots, x_{2 i_2}) 
      \cap (x_{41}, \ldots, x_{4 i_4}). 
  \end{aligned}
\end{displaymath}
Hence, $I$ is the sum of the two squarefree monomial ideals 
$I_1:=(x_{11}, \ldots, x_{1 i_1}) \cap (x_{31}, \ldots, x_{3 i_3})$ and
$I_2:=(x_{21}, \ldots, x_{2 i_2}) \cap (x_{41}, \ldots, x_{4 i_4})$ with 
$\arithdeg I_i = \indeg I_i$ ($i=1, 2$). 
It is known by Schmitt and Vogel \cite[Theorem 1, p.\  247]{SchmVo} that 
$\ara I_i = \pd_S S/I_i$ ($i=1, 2$). 
Since $X (I_1) \cap X (I_2) = \emptyset$, 
we have $\ara I = \pd_S S/I$ by Proposition \ref{claim:connected}. 

\par
For the case of ${\calH}_{14}$, 
we use Lemma \ref{claim:SchmVo} again. 
In this case, $\pd_S S/I = N-2$. 
Note that $N = i_1 + \cdots + i_6$ since $j_2 = j_3 = j_4 = 0$. 
We define 
${\mathcal{P}}_{\ell} := {\mathcal{P}}_{\ell}^{(1)} \cup \cdots \cup 
{\mathcal{P}}_{\ell}^{(5)}$, $\ell = 0, 1, \ldots, N-3$, 
by 
\begin{displaymath}
  \begin{aligned}
    {\mathcal{P}}_{\ell}^{(1)} 
    &= \left\{ x_{1 {\ell}_1} x_{3 {\ell}_3} \; : \; 
      {\ell}_1 + {\ell}_3 = \ell + 2; \  
      1 \leq {\ell}_s \leq i_s \  (s = 1, 3)
    \right\}, \\
    {\mathcal{P}}_{\ell}^{(2)} 
    &= \left\{ x_{1 {\ell}_1} w_{3 {\ell}_3} w_{4 {\ell}_4} \; : \; 
    \begin{aligned}
      &{\ell}_1 + {\ell}_3 + {\ell}_4 + i_3 = \ell + 3, \\
      &1 \leq {\ell}_1 \leq i_1; \  
       1 \leq {\ell}_3 \leq i_2 + i_5; \  
       1 \leq {\ell}_4 \leq i_4 + i_6, \\
      &\text{${\ell}_3 \leq i_2$ or ${\ell}_4 \leq i_4$}  
    \end{aligned}
    \right\}, \\
    {\mathcal{P}}_{\ell}^{(3)} 
    &= \left\{ x_{3 {\ell}_3} w_{1 {\ell}_1} w_{2 {\ell}_2} \; : \; 
    \begin{aligned}
      &{\ell}_3 + {\ell}_1 + {\ell}_2 + i_1 = \ell + 3, \\
      &1 \leq {\ell}_3 \leq i_3; \  
       1 \leq {\ell}_1 \leq i_4 + i_5; \ 
       1 \leq {\ell}_2 \leq i_2 + i_6, \\
      &\text{${\ell}_1 \leq i_4$ or ${\ell}_2 \leq i_2$}  
    \end{aligned}
    \right\}, \\
    {\mathcal{P}}_{\ell}^{(4)} 
    &= \left\{ x_{2 {\ell}_2} x_{4 {\ell}_4} \; : \; 
      {\ell}_2 + {\ell}_4 + i_1 + i_3 = \ell + 2; \  
      1 \leq {\ell}_s \leq i_s \  (s = 2, 4)
    \right\}, \\
    {\mathcal{P}}_{\ell}^{(5)} 
    &= \left\{ x_{5 {\ell}_5} x_{6 {\ell}_6} \; : \; 
      {\ell}_5 + {\ell}_6 + i_1 + \cdots + i_4 = \ell + 3; \  
      1 \leq {\ell}_s \leq i_s \  (s = 5, 6)
    \right\}. \\
  \end{aligned}
\end{displaymath}
Here, 
\begin{displaymath}
  \begin{alignedat}{3}
    w_{1 {\ell}_1} &= \left\{ 
    \begin{alignedat}{3}
      &x_{4 {\ell}_1}, &\quad & 1 \leq {\ell}_1 \leq i_4, \\
      &x_{5 {\ell}_1 - i_4}, &\quad & i_4 < {\ell}_1 \leq i_4 + i_5, 
    \end{alignedat}
    \right. 
    &\qquad 
    w_{2 {\ell}_2} &= \left\{ 
    \begin{alignedat}{3}
      &x_{2 {\ell}_2}, &\quad & 1 \leq {\ell}_2 \leq i_2, \\
      &x_{6 {\ell}_2 - i_2}, &\quad & i_2 < {\ell}_2 \leq i_2 + i_6, 
    \end{alignedat}
    \right. \\
    w_{3 {\ell}_3} &= \left\{ 
    \begin{alignedat}{3}
      &x_{2 {\ell}_3}, &\quad & 1 \leq {\ell}_3 \leq i_2, \\
      &x_{5 {\ell}_3 - i_2}, &\quad & i_2 < {\ell}_3 \leq i_2 + i_5, 
    \end{alignedat}
    \right. 
    &\qquad 
    w_{4 {\ell}_4} &= \left\{ 
    \begin{alignedat}{3}
      &x_{4 {\ell}_4}, &\quad & 1 \leq {\ell}_4 \leq i_4, \\
      &x_{6 {\ell}_4 - i_4}, &\quad & i_4 < {\ell}_4 \leq i_4 + i_6. 
    \end{alignedat}
    \right. 
  \end{alignedat}
\end{displaymath}

\par
Note that for each $k$ $(k=1, 2, \ldots, 5)$, 
the range of $\ell$ with
${\mathcal{P}}_{\ell}^{(k)} \neq \emptyset$ is given by the 
following list: 
\begin{displaymath}
  \begin{aligned}
    {\mathcal{P}}_{\ell}^{(1)}: \quad &0 \leq \ell \leq i_1 + i_3 -2, \\
    {\mathcal{P}}_{\ell}^{(2)}: \quad 
       &i_3 \leq \ell \leq \max \{ N - i_5 - 3, N - i_6 - 3 \}, \\
    {\mathcal{P}}_{\ell}^{(3)}: \quad 
       &i_1 \leq \ell \leq \max \{ N - i_5 - 3, N - i_6 - 3 \}, \\
    {\mathcal{P}}_{\ell}^{(4)}: \quad 
       &i_1 + i_3 \leq \ell \leq i_1 + \cdots + i_4 -2, \\
    {\mathcal{P}}_{\ell}^{(5)}: \quad 
       &i_1 + \cdots + i_4 - 1\leq \ell \leq i_1 + \cdots + i_6 -3 = N - 3. 
  \end{aligned}
\end{displaymath}

\par
By a similar way to the case of ${\calH}_{17}$, we can check that 
${\mathcal{P}}_{\ell}$, $\ell = 0, 1, \ldots, N-3$ satisfy 
the conditions (SV1), (SV2), and (SV3). 
For the condition (SV3), as in ${\calH}_{17}$, 
we list the choices for $a' \in {\mathcal{P}}_{{\ell}'}^{(k')}$ 
from given $a, a'' \in {\mathcal{P}}_{\ell}$ in Table \ref{tab:H3_list}. 

\begin{table}
  \caption{The choices for $a'$ in the case of ${\calH}_{14}$.}
  \label{tab:H3_list}
  \begin{center}
    \begin{tabular}{ccccl} \hline
      $k$ & $k''$ & $k'$ & additional condition & $a'$ \\ \hline
      $1$ & $2$   & $1$  & & $x_{1 {\ell}_1''} x_{3 {\ell}_3}$ \\
      $1$ & $3$   & $1$  & & $x_{1 {\ell}_1} x_{3 {\ell}_3''}$ \\
      $1$ & $4,5$ & \multicolumn{3}{c}{These cases do not occur.} \\
      $2$ & $3$   & $1$  & & $x_{1 {\ell}_1} x_{3 {\ell}_3''}$ \\
      $2$ & $4$   & $2$  & & 
        $x_{1 {\ell}_1} x_{2 {\ell}_2''} x_{4 {\ell}_4''} 
         = x_{1 {\ell}_1} w_{3 {\ell}_2''} w_{4 {\ell}_4''}$ \\
      $2$ & $5$   & $2$  & ${\ell}_3 \leq i_2$ & 
        $x_{1 {\ell}_1} w_{3 {\ell}_3} x_{6 {\ell}_6''} 
         = x_{1 {\ell}_1} w_{3 {\ell}_3} w_{4 i_4 + {\ell}_6''}$ \\
      $2$ & $5$   & $2$  & ${\ell}_4 \leq i_4$ & 
        $x_{1 {\ell}_1} x_{5 {\ell}_5''} w_{4 {\ell}_4} 
         = x_{1 {\ell}_1} w_{3 i_2 + {\ell}_5''} w_{4 {\ell}_4}$ \\
      $3$ & $4$   & $3$  & & 
        $x_{3 {\ell}_3} x_{4 {\ell}_4''} x_{2 {\ell}_2''} 
         = x_{3 {\ell}_3} w_{1 {\ell}_4''} w_{2 {\ell}_2''}$ \\
      $3$ & $5$   & $3$  & ${\ell}_1 \leq i_4$ & 
        $x_{3 {\ell}_3} w_{1 {\ell}_1} x_{6 {\ell}_6''} 
         = x_{3 {\ell}_3} w_{1 {\ell}_1} w_{2 i_2 + {\ell}_6''}$ \\
      $3$ & $5$   & $3$  & ${\ell}_2 \leq i_2$ & 
        $x_{3 {\ell}_3} x_{5 {\ell}_5''} w_{2 {\ell}_2} 
         = x_{3 {\ell}_3} w_{1 i_4 + {\ell}_5''} w_{2 {\ell}_2}$ \\
      $4$ & $5$ & \multicolumn{3}{c}{This case does not occur.} \\ \hline
    \end{tabular}
  \end{center}
\end{table}

\par
This completes the proof of Theorem \ref{claim:Adual_4gen}. 

\begin{remark}
  \label{rmk:H1toH3}
  For the exceptional case $\calH = \calH_{14}$, 
  we can also consider the corresponding elements on the construction 
  for $\calH_{17}$; in fact, these generate $I$ up to radical. 
  But in this case, the range of $\ell$ with 
  ${\mathcal{P}}_{\ell} \neq \emptyset$ is 
  $0 \leq \ell \leq N-2$ because of ${\mathcal{P}}_{\ell}^{(7)}$. 
  Therefore, we cannot obtain $\ara I$ from this construction. 
  This is also true for $\calH = \calH_1$. 
\end{remark}

\begin{acknowledgement}
  The authors thank the referee for reading the manuscript carefully. 
\end{acknowledgement}


\end{document}